\documentclass[11pt,leqno]{article}
\usepackage{amsthm,amsfonts,amssymb,epsfig,graphics,amsmath,eufrak}
 \usepackage[latin1]{inputenc}\relax

\newcommand{\tyek}{t^{-1}}
\newcommand{\tnim}{t^{-\frac 12}}
\newcommand{\trob}{t^{-\frac 14}}

\newcommand{\tpnim}{t^{-\frac 12 (1-\frac 1p)-\frac 12}}

\newcommand{\ubar}{\bar{u}}
\newcommand{\utild}{\tilde{u}}
\newcommand{\fe}{\varphi}
\newcommand{\hcal}{\mathcal{H}}
\newcommand{\fcal}{\mathcal{F}}
\newcommand{\gtild}{\tilde{G}}
\newcommand{\pstar}{p^*}
\newcommand{\dcal}{\mathcal{D}}
\newcommand{\mcal}{\mathcal{M}}
\newcommand{\ddp}{\frac{\partial\ubar^\delta}{\partial
\delta}}
\newcommand{\ddpu}{\frac{\partial\bar U^\delta}{\partial
\delta}}
\newcommand{\ddpw}{\frac{\partial\bar W^\delta}{\partial
\delta}}
\newcommand{\ddpwx}{\frac{\partial\bar W_x^\delta}{\partial
\delta}}

\newcommand{\feh}{\hat\fe}
\newcommand{\BbbR}{\Bbb{R}}
\newcommand{\CalO}{\mathcal{O}}
\newcommand{\ta}{\tilde{A}}
\newcommand{\tao}{\tilde{A}^0}
\newcommand{\tw}{\tilde{W}}
\newcommand{\tb}{\tilde{B}}
\newcommand{\ba}{\bar{A}}

\newcommand{\bb}{\bar{B}}
\newcommand{\bw}{\bar{W}}
\newcommand{\tu}{\tilde{U}}
\newcommand{\bu}{\bar{U}}
\newcommand{\duav}{dU_{Ave}}
\newcommand{\emi}{\mathcal{M}_1}
\newcommand{\emii}{\mathcal{M}_2}

\def\D{\partial }
\newcommand\adots{\mathinner{\mkern2mu\raise1pt\hbox{.}
\mkern3mu\raise4pt\hbox{.}\mkern1mu\raise7pt\hbox{.}}}

\newtheorem{theo}{Theorem}[section]
\newtheorem{prop}[theo]{Proposition}
\newtheorem{cor}[theo]{Corollary}
\newtheorem{lem}[theo]{Lemma}

\newtheorem{exam}[theo]{Example}
\newtheorem{rem}[theo]{Remark}

\numberwithin{equation}{section}

\title{
\bf $L^p$ Asymptotic Behavior of Perturbed Viscous Shock Profiles
}

\author{\sc \small
 \sc \small Mohammadreza Raoofi\thanks{Indiana University, Bloomington, IN 47405;
mraoofi@indiana.edu. This work was carried out as part of the
author's doctoral thesis at Indiana University, Bloomington, under
the direction of Kevin Zumbrun. The author would like to thank
professor Zumbrun for his advice, encouragement and support. This
project was supported in part by National Science Foundation under
Grant DMS-0070765 }
 }

\begin{document}

\maketitle

\begin{abstract}
We investigate the $L^p $ asymptotic behavior $( 1\le p \le
\infty)$ of a perturbation of a Lax or overcompressive type shock
wave solution to a system of conservation law in one dimension.
 The system of
the equations can be strictly parabolic, or have real viscosity
matrix (partially parabolic, e.g., compressible Navier--Stokes
equations or equations of Magnetohydrodynamics). We use  known
pointwise Green function bounds for the linearized equation around
the shock
 to show that the perturbation of such a
solution can be decomposed into 
a part corresponding to shift in shock position or shape, 
a part which is the sum of
diffusion waves, i.e., the solutions to a viscous Burger's
equation, conserving the initial mass and convecting away from the
shock profile in outgoing modes, 
and another part which is more rapidly decaying in $L^p$.
\end{abstract}



\section{Introduction} \label{S:introduction}
Consider the system of conservation laws with viscosity:
\begin{equation}
u_t + f(u)_x =\nu (B(u)u_x)_x \label{conservationlaw}
\end{equation}
with $u\in \mathbb{R}^n$ is the  conserved quantity, and $\nu$ is
a constant measuring transport effects (e.g. viscosity or heat
conduction). As we are not considering the vanishing-viscosity
limit $\nu \to 0$, we can assume $\nu=1.$
An important class of the solutions for (\ref{conservationlaw})
are the viscous shock wave solutions, i.e., solutions in the form
$\ubar(x,t) = \ubar(x-st)$, where the constant $s$ is the velocity
of the shock, and where $\ubar$ connects the endstates
$u_\pm=\ubar(\pm\infty)$. With a simple change of coordinates, we
can assume that $s=0$ (a \textit{stationary} shock solution).
$\ubar$ is assumed to be an element of a smooth manifold
$\{\ubar^\delta\}_{\delta\in\mathbb{R}^d}$, which consists of
stationary solutions of (\ref{conservationlaw}) connecting the
same endstates $u_-$ and $u_+$, and $\ubar=\ubar^0$. Loosely
stated, we prove that, with $\utild$ a solution of
(\ref{conservationlaw}) and a small perturbation of $\ubar,$ there
is a small $\delta_0$ such that
\begin{equation}
\utild(x,t)-\ubar^{\delta_0}(x)= v(x,t)+\fe(x,t)+\ddp\delta(t)
\label{decompose}
\end{equation}
where
\begin{enumerate}
\item $ |v(\cdot, t)|_{L^p} \sim
(1+t)^{-\frac12(1-\frac1p)-\frac14}$. By choosing appropriate
$\delta_0$ and the mass carried by $\fe$, we also obtain  zero
initial mass for $v$, i.e.,\linebreak
$\int_{-\infty}^{+\infty}v(x,0)\,dx=0.$
\item $\fe$ is a summation of convecting diffusion waves,
i.e., self similar solutions to the viscous Burger's equation with
appropriate coefficients, propagating away from the shock and
preserving the initial mass in the outgoing modes, $|\fe(\cdot,
t)|_{L^p}\sim (1+t)^{-\frac12(1-\frac1p)}$,
 and
\item$|\delta(t)|\sim (1+t)^{-\frac12+\epsilon}$ and
$\delta(0)=0$. One can view $\delta(t)$  as indexing the ``instantaneous"
shock location and shape: employing Taylor's expansion gives us
\begin{equation}
\ubar^{\delta_0+\delta(t)}-\ubar^{\delta_0}=\ddp\delta(t) + \bold
O (|\delta(t)|^2 e^{-k|x|})\label{ishl}
\end{equation}
which shows that $\ddp\delta(t)$ corresponds to a shift in the
shock location and or shape (up to an error of order $|\delta(t)|^2 e^{-k|x|}$,
which  decays faster than $v$ in any $L^p$ norm).
\end{enumerate}
For the exact definitions and conditions, see the subsequent
sections; especially see Theorem \ref{thmmainover}, corollaries
\ref{cor1}, \ref{cor2} and their counterparts in the real
viscosity case, 
which comprise the main results of this paper.

 To prove the above statements, we use (\ref{decompose})
and initial equations for $\utild$ and $\ubar$ to obtain
\begin{equation}
v_t-Lv=\mathcal{R} (v, v_x, \ddp\delta(t), \fe)
\end{equation}
where $L$ stands for the linearized operator around $\ubar$ and
$\mathcal{R}$ is a remainder we get applying Taylor's expansion.
If $G(x,t;y)$ is the Green function corresponding to ${\D_t} - L$,
then applying Duhamel's principle yields:
\begin{equation}
\begin{aligned}
v(x, t)=&\int_{-\infty}^{+\infty}G(x,t;y)v(y,0)dy\\
&+\int_0^t\int_{-\infty}^{+\infty}G(x,t-s;y)\mathcal{R}(v,
\ddp\delta(t), \fe)(y,s)dyds;
\end{aligned} \label{sabz}
\end{equation}
we try, then, to use a continuous induction to prove the desired
rates of decay for $v$.

The observation that a perturbation of a shock wave solution to
(\ref{conservationlaw}) can be decomposed into a sum of diffusion
waves and a more rapidly decaying term is due to T.P. Liu (see
\textbf{\cite{Liu,Liu85}}). He proved the result for  weak shocks
and with the viscosity matrix $B$ being identity matrix. To prove
the result, he first constructed an approximate Green function
using heavily the weak shock wave assumption and the identity
matrix, and then used an elaborate pointwise nonlinear iteration
scheme.

We, on the other hand, have already at our disposal the Green
function bounds we need. These sharp bounds are the result of a
 ``dynamical system" approach based on Evans function and
inverse Laplace transform techniques. This approach began for the
strictly parabolic case (\textbf{\cite{GZ, ZH}}) and then was
extended to many other, more physical, regimes, such as real
viscosity case (\textbf{\cite{MaZ.1}} -- \textbf{\cite{MaZ.4}};
see also \textbf{\cite{Z, Z.4})}. In these papers the $L^p, 1<p\le
\infty$, asymptotic stability and $L^1$ asymptotic boundedness, of
Lax type shock profiles  were stated and proved by finding sharp
Green function bounds for the linearized equation; some hints have
also been given about Green function bounds for the
overcompressive case. This approach does not require any
assumption of the weakness of the shock profile, and the
structural and technical assumptions made about the equation and
the wave are rather minimal. However, less information than what
Liu's approach yields has been obtained about the behavior of the
perturbation.

Using the very same Green distribution bounds, we prove the
results which Liu first observed with fewer assumptions: the shock
profile can have small or large amplitude, be of Lax or
overcompressive type (to our knowledge, this result is the first
rather complete result about asymptotic behavior of a perturbed
overcompressive shock); the viscosity term $(B(u)u_x)_x$ can be
strictly hyperbolic
 (section \ref{S:overc}) or have the block structure of the
real viscosity case (section \ref{S:realvisc}). Also no pointwise
bounds on initial data 
are 
required, only bounds on $L^p$ norm and
moments. In return for localization of the initial data, Liu
obtains pointwise bounds for the perturbation. We only assume
smallness of the initial perturbation and its moment in some $L^P$
spaces, but then no pointwise information is obtained (we believe,
however, that with a similar method and some more work, and with
localization of the initial data, we can achieve pointwise bounds
similar to those Liu has obtained).

\textbf{Plan of the paper:}  In computing the bounds for $v$ using
(\ref{sabz}), we frequently use Young--Hausdorf inequality and
$L^p$ norms of the different components of $\mathcal{R}$. However,
a term that does not yield the necessary estimates this way occurs
and that is when part of the Green function that is like a
convecting heat kernel is convoluted against $(\fe_i^2)_y(y,s)$,
with $\fe_i$ a diffusion wave convecting at a different speed.
Sharp estimation of such terms was first treated by Liu {\bf \cite{Liu}}.
Here, we extract the essential features of his argument,
to establish that similar bounds hold whenever
the derivative of $G$ along characteristic directions
decays more rapidly than $G_y$: in particular, for the
more general Green function terms we consider here.
Most of section \ref{S:preliminary} is to compute pointwise bounds
for this part of the calculations. In section \ref{S:constant} we
consider the  case, already well established (see
\textbf{\cite{LZe, Kaw86, Kaw87, CL}}), of a strictly parabolic
system with a perturbation of a constant state solution. The
calculations foreshadow the more difficult case of shock waves. In
sections \ref{S:overc} and \ref{S:realvisc} we consider the
perturbation of a shock wave solution in the strictly parabolic
case and the real viscosity case, respectively, which are very
similar; the main difference is that while in section
\ref{S:overc} we use strict parabolicity to establish short time
estimates and thereby find good bounds for $v_x$, in the real
viscosity case we have to go through a long list of energy
estimates to find the bounds we need for the derivatives.
These bounds generalize similar energy estimates obtained
in {\bf \cite{MaZ.2, MaZ.4, Z.4, Z.5}}, which in turn generalize
the important estimates obtained by Kawashima and others
(see {\bf \cite{Kaw}} and references therein)
for perturbations of constant states.

\section{Some preliminary computations} \label{S:preliminary}
In this section, we establish some pointwise bounds for the
solution $u(x,t)$ of
\begin{equation}\label{avali}
\begin{cases}
u_t - u_{xx}=(K(x-t, t)^2)_x \qquad    &\text{for} \,\,t>0 \\
u(x,0)\,\,=0
\end{cases}
\end{equation}
 Here and elsewhere in this article $K(x,t)$ and $g(x, t)$
 both denote the heat kernel: $g(x,t)= K(x,t)=(4\pi t)^{-\frac
12}e^\frac{-x^2}{4t}$

Using Duhamel's principle we obtain from (\ref{avali}),

\begin{equation}
\begin{aligned} \label{dovomi}
u(x,t)&= \int_0^t \int_{-\infty}^{+\infty}g(x-y, t-s)(K(y-s,
s)^2)_y \, dy \, ds\\
&= \int_0^t \int_{-\infty}^{+\infty}g_y(x-y, t-s) K(y-s, s)^2 \,
dy \, ds
\end{aligned}
\end{equation}
The following bounds for $u$ are essential for obtaining $L^1$
bounds in subsequent sections. They are similar in nature to the
bounds given by \textbf{\cite{Liu}}, but the proof given here is
different from that of Liu, and is somewhat more general.
\begin{prop} \label{sect1main}
Let $u(x,t)$ be the solution of (\ref{avali}) given by
(\ref{dovomi}). If $t \geq 1$, then we have
\begin{equation}
\begin{split}
|u(x,t)| \leq C \trob \left( g(x, 4t) +
g(x-t, 4t) \right)\\
+ C
\chi_{\{\sqrt{t} \leq x \leq t-\sqrt{t}\}}
( t^{-1}x^{- \frac 12} + t^{-\frac 12} (t-x)^{-1}).
\end{split}\label{liubounds}
\end{equation}
 where $\chi$ stands for the
indicator function, and $C$ is a constant independent of $t$
and $x$.\\
The same result holds if $K^2$ in (\ref{avali}) and (\ref{dovomi})
is replaced with $K_x$.
\end{prop}
\begin{rem}\textup{
The result just mentioned is achieved by detecting crucial
cancelation in the calculations. In fact, if we replace the
integrands in (\ref{dovomi}) by their absolute value, we will
obtain the following bounds instead (See \textbf{\cite{HZ}}):
$$|u(x,t)| \leq C  \left( g(x, 4t) +
g(x-t, 4t) \right)\\
+ C \chi_{\{\sqrt{t} \leq x \leq t-\sqrt{t}\}} (  x^{-\frac 12}
(t-x)^{-\frac12}).$$
 }
\end{rem}
\begin{proof}[Proof of Proposition \ref{sect1main}]
As $K^2 \sim K_x$, the proof will be stated only for $K^2$. It
would be straightforward to observe that the same argument works
for $K_x$ at every step.

 We begin the proof by first stating a simple
lemma:
\begin{lem} \label{stlem}
If \, $0\leq s\leq \sqrt{t}$, \, then $e^{\frac{-(x\pm s)^2}{4t}}
\leq Ce^{\frac{-x^2}{8t}}$ with $C$ independent of $t,s$ and $x$.
\end{lem}
\begin{proof}
The statement of the lemma is equivalent to
\begin{equation}
\frac{-(x\pm s)^2}{4t} \leq \frac{-x^2}{8t}+D \notag
\end{equation}
for some $D$, which (after some calculation) in its turn is
equivalent to \linebreak $(x\pm 2s)^2-2s^2 \geq -8Dt$, which holds
for $D > \frac 14$, since $s^2 < t$.
\end{proof}
In the proof of the proposition, we will, in several places, use
the following lemma, which is due to P. Howard \textbf{\cite{Ho}}.
\begin{lem}[P. Howard] \label{holem}
Let $f(\sigma) \geq 0$ be a nonincreasing function on
$\mathbb{R}_{+}$ and $f(0) < \infty$. Assume further that there
exist constants $\gamma > 0, \omega > 1 $  so that $f(\sigma) \geq
\gamma e^{-\frac a2 (1-\frac 1{\omega})^2 \sigma^2}$ on
$\mathbb{R}_{+}$. Then for $a, z > 0$,
$$\int_{0}^{+\infty} e^{-a (z-\sigma)^2} f(\sigma)\, d\sigma \leq \frac{C(\omega)}{\sqrt{a}}
f\Big(\frac{z} {\omega}\Big).$$
\end{lem}
\begin{proof}
See page 102 of \textbf{\cite{Ho}}.
\end{proof}
\begin{rem}\textup{
This result is sharp down to scale $a^{-\frac12}$. In practice,
this is often augmented with $L^{\infty}$ bounds obtained by other
means. See for example lemma 5 of  \textbf{\cite{HZ}} or remark
\ref{liuinitial} of this paper.}
\end{rem}
 In what follows we will frequently use the following properties of the
heat kernel $g$, which are easy to prove:
\begin{align}
 \int_{-\infty}^{+\infty}g(x-y, t)g(y,
t') dy = g(x,t+t') \label{gg1}\\
 |g_x(x,t)|\leq C\tnim g(x,2t)
\label{gg2}\\
 |g_t(x,t)|\leq C\tyek
g(x,2t) \label{gg3}\\ |g(x,t)| \leq C\, \tnim. \label{gg4}
\end{align}
Rewriting (\ref{dovomi}), we have:
\begin{equation}
\begin{aligned}
u(x,t)&= \int_0^t \int_{-\infty}^{+\infty}g(x-y, t-s)(g(y-s,s)^2)_y \, dy \, ds\\
&=\int_0^{\sqrt{t}} \int_{-\infty}^{+\infty}g(x-y, t-s)(g(y-s,
s)^2)_y \, dy \, ds\\
&+ \int_{\sqrt{t}}^{t-\sqrt{t}} \int_{-\infty}^{+\infty}g(x-y,
t-s)(g(y-s, s)^2)_y \, dy \, ds\\
&+ \int_{t-\sqrt{t}}^t \int_{-\infty}^{+\infty}g(x-y, t-s)(g(y-s,
s)^2)_y \, dy \, ds\\
&=: I + II + III.
\end{aligned} \label{shekaf}
\end{equation}
$(I)$ and $(III)$ are easy to estimate:
\begin{align}
|I| &= \left| \int_0^{\sqrt{t}} \int_{-\infty}^{+\infty}g(x-y,
t-s)(g(y-s, s)^2)_y \, dy \, ds\right| \label{I1}\\
&= \left| \int_0^{\sqrt{t}} \int_{-\infty}^{+\infty}g_y(x-y, t-s)
g(y-s, s)^2 \, dy \, ds\right|.\label{I2}
\end{align}
By (\ref{gg2}) and (\ref{gg4}), the above is,
$$ \leq  C\int_0^{\sqrt{t}}
\int_{-\infty}^{+\infty}(t-s)^{-\frac12} s^{-\frac12} g(x-y,
2(t-s))\,g(y-s, 2s) \, dy \, ds $$
which is, by (\ref{gg1}),
$$
 \leq  C\int_0^{\sqrt{t}}
(t-s)^{-\frac12} s^{-\frac12} g(x-s, 2t) \, ds. \notag
$$
Now, using lemma \ref{stlem}, the above is
\begin{eqnarray}
 &\leq&  C\,g(x, 4t) \int_0^{\sqrt{t}} \notag
(t-s)^{-\frac12} s^{-\frac12} \, ds\\ \notag &\leq&  C\tnim g(x,
4t) \int_0^{\sqrt{t}} s^{-\frac12} \, ds\\ \notag
 &\leq& C\trob g(x,4t).
\end{eqnarray}
 Part $(III)$
in (\ref{shekaf}) can be handled similarly.

The more difficult part is part $(II)$ of (\ref{shekaf}):
\begin{equation}
II= \int_{\sqrt{t}}^{t-\sqrt{t}} \int_{-\infty}^{+\infty}g(x-y,
t-s)(g(y-s, s)^2)_y \, dy \, ds.
\end{equation}
In order to estimate $(II)$, let us write $g=g(x, \tau)$. We have
then,
\begin{equation}
\begin{aligned}
g(x-y, t-s)(g(y-s, s)^2)_y \,&= (g(x-y, t-s)g(y-s, s)^2)_s \\ & -
\,g_\tau (x-y, t-s)g(y-s, s)^2\\ &+\, g(x-y, t-s)(g^2)_\tau(y-s,
s) .
\end{aligned} \label{shekaf2}
\end{equation}
We will do the estimates piece by piece.

The first part of (\ref{shekaf2}) can be estimated as follows.
\begin{align}
\int_{\sqrt{t}}^{t-\sqrt{t}} \int_{-\infty}^{+\infty} (g(x-y,
t-s)g(y-s, s)^2)_s\, dy\, ds\\
= \int_{-\infty}^{+\infty} g(x-y, \sqrt{t})\,g(y-t+\sqrt{t},
t-\sqrt{t})^2\,dy\\
- \int_{-\infty}^{+\infty} g(x-y, t-\sqrt{t})\,g(y-\sqrt{t},
\sqrt{t})^2\,dy.
\end{align}
Using (\ref{gg1}) and (\ref{gg4}), we will have:
$$ \int_{-\infty}^{+\infty} g(x-y, \sqrt{t})\,g(y-t+\sqrt{t},
t-\sqrt{t})^2\,dy \leq C\tnim g(x-t+\sqrt{t}, t), $$ and
$$\int_{-\infty}^{+\infty} g(x-y, t-\sqrt{t})\,g(y-\sqrt{t},
\sqrt{t})^2\,dy \leq C\trob g(x-\sqrt{t}, t),$$ but by lemma
\ref{stlem}
$$g(x-\sqrt{t}, t) \leq g(x, 2t)$$
$$g(x-t+\sqrt{t}, t) \leq g(x-t, 2t).
$$These terms fit in the right hand side of (\ref{liubounds}).

For the other parts in (\ref{shekaf2}), we use (\ref{gg3}) to
obtain:
\begin{equation}
\begin{aligned}
&\left| \int_{\sqrt{t}}^{t-\sqrt{t}} \int_{-\infty}^{+\infty}
g_\tau (x-y, t-s) g(y-s, s)^2 \, dy\, ds
\right|  \\
\leq &\int_{\sqrt{t}}^{t-\sqrt{t}} s^{-\frac 12} (t-s)^{-1} g(x-s,
2t) ds  \\
=&\int_{\sqrt{t}}^{\frac t2} s^{-\frac 12} (t-s)^{-1} g(x-s, 2t)
ds \\
+ &\, \int_{\frac t2}^{t-\sqrt{t}} s^{-\frac 12} (t-s)^{-1}
g(x-s, 2t) ds \\
=: \,&\mathcal{A} \,+\, \mathcal{B}.
\end{aligned} \label{AB}
\end{equation}

If $x \leq \sqrt{t}$, then
\begin{equation}
\begin{aligned}
\mathcal{A} &\leq C\,\tyek g(x-\sqrt{t}, 2t)\int_{\sqrt{t}}^{\frac t2} s^{-\frac 12}\, ds\\
&\leq  C \, t^{-\frac 12} \, g(x-\sqrt{t}, 2t)\\
&\leq C \tnim g(x, 4t).
\end{aligned}
\end{equation}
Similarly when $x\ge t-\sqrt{t}$, we obtain $\mathcal A \le C\tnim
g(x-t, 4t).$

Now for $\sqrt{t} \leq x \leq \frac t2$ we have:
\begin{align}
\mathcal{A} &= \int_{\sqrt{t}}^{\frac t2} s^{-\frac 12}
(t-s)^{-1} g(x-s, 2t) ds \\
&\leq t^{-1} \int_{\sqrt{t}}^{\frac t2} s^{-\frac 12}
 g(x-s, 2t) ds \\
&= C\, \tyek \int_{\sqrt{t}}^{\frac t2} \Big(\frac st\Big)^{-\frac 12}
 e^{-\frac t8 (\frac xt -\frac st)^2}\, \frac
{ds}t\\
&= C\, \tyek \int_{\frac{1}{\sqrt{t}}}^{\frac 12} {\sigma}^{-\frac
12}
 e^{-\frac t8 (\frac xt -\sigma)^2}\, d\sigma
\label{iin}
\end{align}
Now we use Howard's lemma, lemma \ref{holem}, to find  that
(\ref{iin}) is indeed,
\begin{align}
& \leq C\, t^{-\frac 32} (\frac xt)^{-\frac 12} \\
&=C \, t^{-1}x^{-\frac12}.
\end{align}

%
If, on the other hand, $\frac t2\le x \le t-\sqrt t$, then clearly
$\mathcal{A}$ is majorized by the value already computed for
$x=\frac t2$, of
\begin{equation}
\mathcal{A} \leq C t^{-\frac32}\le C \, t^{-1}x^{-\frac12},
\end{equation}
also acceptable.

Part $\mathcal{B} $ in (\ref{AB}) can be estimated similarly.

Now remains the last part of (\ref{shekaf2}), i.e.,
$$\left| \int_{\sqrt{t}}^{t-\sqrt{t}}\int_{-\infty}^{+\infty} g(x-y, t-s)(g^2)_\tau(y-s,
s) \, dy\, ds\right|,$$ which can easily  be shown to be
\begin{equation}
\leq \int_{\sqrt{t}}^{t-\sqrt{t}}s^{-\frac 32} g(x-s, 2t)\, ds.
\label{aan}
\end{equation}

If $x\,\leq\, \sqrt{t}$, then (\ref{aan}) is
$$\leq C\, g(x-\sqrt{t}, 2t)\int_{\sqrt{t}}^{t-\sqrt{t}}s^{-\frac 32}\, ds$$
$$\leq C\,\trob \, g(x-\sqrt{t}, 2t)$$
$$\leq C\,\trob \, g(x,4t). $$
For $x\, \geq \, t-\sqrt{t}$ we use a similar method.\\

For $\sqrt{t} \, \leq \,x\, \leq \, t-\sqrt{t}$, we use a similar
method to what we used for the previous case, and again invoke the
result in lemma \ref{holem} to conclude that $(\ref{aan}) \, \leq
\,
x^{-\frac 32}$\\
 This completes our proof.
\end{proof}
As an straightforward consequence to the above proposition we
have:
\begin{cor} \label{sect1cor}
For $u(x,t)$ the solution of (\ref{avali}) given by
(\ref{dovomi}), and for $t \geq 1$, we have:
\begin{equation}
\begin{aligned}
|u(\cdot,t)|_{L^p} \leq Ct^{-\frac 12 (1-\frac 1p)-\frac 14},
\end{aligned}
\end{equation}
where $|\cdot|_{L^p}$ stands for the norm in $L^p(\mathbb{R})$,
$1\leq
p\leq +\infty$, and $C$ is a constant.\\
The same result holds if $K^2$ in (\ref{avali}) and (\ref{dovomi})
is replaced with $K_x$.
\end{cor}

Bounds obtained in corollary (\ref{sect1cor}) will be used when
working to prove $L^p$ bounds in the system case. In that setting,
however, we will usually have a  convecting diffusion wave instead
of $K(x , t)$ in the source, and part of a Green function,
convecting at a different speed, in place of $g(x,t)$ in
(\ref{dovomi}). The following result deals with those cases.

\begin{cor} \label{phicor}
In (\ref{dovomi}), if we replace $K(y-at, t)$ with $\phi(y,t)$,
and $g(x-y, t)$ with $G(x,t;y)$ in (\ref{dovomi}), we will have
similar bounds for $u(x,t)$ obtained in proposition
(\ref{sect1main}) and corollary (\ref{sect1cor}), provided $\phi$
and $G$ satisfy the following bounds:
\begin{align}
|G(x,t;y)| \leq C g(x-y-at, \beta t),\\
|G_y(x,t;y)| \leq C \tnim g(x-y-at, 2\beta t),\\
|G_t(x,t;y)| \leq C \tyek g(x-y-at, 2\beta t),\\
|\phi(x,t)| \leq C g(x-bt, \beta t), \label{bd3}\\
|\phi_y(x,t)| \leq C \tnim g(x-bt,2\beta t),\\
|\phi_t(x,t;y)| \leq C \tyek g(x-bt, 2\beta t), \label{bd6}
\end{align}
for some $a\neq b,$ and some constants $C, \beta >0$.
\end{cor}
\begin{proof}
A review of the proofs just carried out will show that the above
bounds are the only ones used in the proof. Hence everything works
in the same way as before.
\end{proof}
\begin{exam} \label{phi}
$$\phi (x,t) = \frac{(e^{m/2\sqrt{\beta}}-1)\tnim e^{\frac{-x^2}{4\beta t}}}
{\frac{2\sqrt\pi}{\sqrt{\beta}} +
(e^{m/2\sqrt{\beta}}-1)\int_{\frac{x}{\sqrt{4\beta
t}}}^{+\infty}e^{-\xi^2}d\xi}$$ \textup{solves}
\begin{equation}
\begin{cases}
\phi_t - \beta \phi_{xx} = -(\phi^2)_x \qquad
\text{for} &t > 0\\
\phi (x, 0) = m\delta_0 &t= 0,
\end{cases}
\end{equation}
\textup{where $\delta_0$ stands for the Dirac distribution, and $m
= \int_{-\infty}^{+\infty} \phi(x, t) \, dx$.\\ The function
$\phi$ satisfies the inequalities in (\ref{phicor}) (see
\textbf{\cite{Liu85}}). Also one can easily see  that if one puts
$\phi_x$ in place of $\phi^2$ in above argument, then again one
will obtain similar results, as $\phi_x \sim \phi^2.$ This
function $\phi$ will be the prototype of the \emph{diffusion
waves}, which we are going to define and use in the next
sections.}
\end{exam}

Now that we are in the mood of working with the heat kernels, let
us state some lemmas that we will need when dealing with systems.
The proofs are easy and left to the reader.

\begin{lem} \label{interaction}
If $a_1 \neq a_2$ and $\beta_1, \beta_2 > 0$, then $|K(x-a_1t,
\beta_1 t) K(x-a_2t, \beta_2t)|_{L^p} \leq C e^{-\eta t},$ for
some
$\eta > 0$.\\
The same result holds if one replaces $K$ with $\phi$ from example
(\ref{phi}).
\end{lem}

And the following lemmas, which will be needed for shock wave
cases:
\begin{lem} \label{hKlem}
Assume $a > 0 $ (respectively, $ a < 0)$, $h(x)$ is a bounded
function and $h(x)= \bold {O} (e^{-|x|})$ as $x \to +\infty$
(respectively,
 as $x \to -\infty).$ Then $|h(x) K(x-at,
t)|_{L^p} = \bold {O} ( e^{-\eta t})$ for some $\eta > 0$.\\
\end{lem}

\begin{lem} \label{KKKlem}
Assume $a, b $ are both of the same sign, then
$$\int_{0}^{t}\int_{-\infty}^{+\infty} g(y + a(t-s), t-s) g(y-bs, s)^2 dy\, ds
=\bold{O} (e^{-\eta t})$$ for some $\eta > 0$; If $a$ and $b$ are
of the different sign then,
$$\int_{0}^{t}\int_{-\infty}^{+\infty} g(y + a(t-s), t-s) g(y-bs, s)^2 dy\, ds
=\bold{O} (t^{-\frac12}).$$
\end{lem}

\begin{lem} \label{Kelem}
 For $t>0$ and for $a \neq 0$,
$$\int_{-\infty}^{+\infty}g(y+at, t)e^{(-|y|)}dy=\bold O (t^{-\frac12}e^{-\eta t}).$$
\end{lem}

\section{System of conservation laws with constant background solution} \label{S:constant}
 Now consider the system of conservation laws
\begin{equation}
\utild_t + f(\utild)_x = (B(\utild) \utild_{x})_x
\end{equation}
with the solution $\utild = \utild (x, t)\in \mathbb{R}^n,$  a
 perturbation of the constant background solution
$\ubar \equiv \ubar_0$.

Let $u = \utild - \ubar$, and use Taylor's expansion to obtain
\begin{equation}
u_t + A u_x - B u_{xx} = - (\Gamma (u, u))_x + \bold {O}(|u|^3)_x
+ \bold {O}(|u||u_x|)_x \label{utassys}
\end{equation}
where $A = df(\ubar)$, $B = B(\ubar)$ and $\Gamma = \frac 12 d^2f(\ubar)$.\\

Some basic assumptions have to be made: we assume $f, B \in C^3$,
$df(\ubar)$ is strictly hyperbolic, $Re \,\sigma(B) > 0$ and
finally stability criterion of Majda and Pego \textbf{\cite{Kaw,
MP}}:
 $Re \, \sigma(-ikA-k^2B)<-\theta k^2$ for all real $k$ and some $\theta
> 0.$\\

Let $a_1, \cdots, a_n$ be the eigenvalues of $A$ (necessarily
disjoint by the strict hyperbolicity of $A)$, with corresponding
right eigenvectors $r_1, \cdots, r_n,$ and left eigenvectors $l_1,
\cdots, l_n$, normalized so that $l_i\cdot r_j = \delta_{ij}.$
Define $\Gamma_{jk}^i$ and $b^i_j$ to be the constant coefficients
satisfying
\begin{equation}
\Gamma (r_j, r_k) = \sum_{i=1}^n \Gamma_{jk}^i r_i,\quad Br_j=
\sum_{i=1}^n b^i_j r_i \label{coefficient}
\end{equation}
 and set $\beta_i = b_i^i$ and $\gamma_i = \Gamma_{ii}^i $ (notice that it
 follows from our assumptions about $A$ and $B$ that $\beta_i > 0)$. Define
 {\it diffusion wave} in $r_i$ direction: $\fe^i( x, t)$ to be the solution
 of
\begin{equation}
\begin{cases}
\fe^i_t + a_i\fe^i_x - \beta_i \fe^i_{xx} = -\gamma_i (\fe^{i2})_x
\qquad
\text{for} &t > -1\\
\fe^i (x, -1) = m_i\delta_0 &t= -1
\end{cases} \label{fiith}
\end{equation}
where $\delta_0$ is the Dirac distribution, and $m_i$ is the
amount of mass $\int_{-\infty}^{+\infty} u(x, 0) \, dx$  in the
direction $r_i$.  See Example (\ref{phi}). Assuming $m\leq E_0$ we
will have:
\begin{equation}
\left|\frac {\partial^n}{\partial x^n}\fe^i(\cdot ,
t)\right|_{L^p} \leq CE_0 (1+t)^{-\frac 12 (1-\frac 1p)-\frac n2},
\label{fibounds}
 \end{equation}
 i.e., $\fe^i$ has $L^p$ bounds like a heat kernel. It is not
difficult to observe that $\fe$ acts like a convecting heat
kernel; especially for our interest is the fact that it satisfies
the bounds (\ref{bd3})-(\ref{bd6}). Finally set $$\fe =
\sum_{i=1}^n \fe^i r_i.$$

 Let $v := u - \fe$, hence
\begin{equation}
\utild = \ubar + \fe + v.
\end{equation}
Notice that
\begin{equation}
 \int_{-\infty}^{+\infty} v_0 dx = 0.  \label{0masssys}
\end{equation}
Set $V_0 (x) := \int_{-\infty}^{x} v_0 dx$.

\begin{lem} \label{wW}
$|V_0|_{L^1} = \int_{-\infty}^{+\infty} |V_0| dx \leq |x v(x,
0)|_{L^1},$ assuming that the latter quantity is bounded.
\end{lem}
\begin{proof}
We can assume, without loss of generality, that $v_0$ is scalar,
i.e., $v_0(x)\in \mathbb{R}.$ If $x_0 < 0$, then $|V_0(x_0)| =
|\int_{-\infty}^{x_0} v_0 (x) dx| \leq \frac 1{x_0}
\int_{-\infty}^{x_0} |x| |v_0 (x)| dx$, hence $|x_0 V_0(x)| \leq
\int_{-\infty}^{x_0} |x| |v_0 (x)| dx$, which approaches $0$, as
$x_0 \rightarrow -\infty$. Likewise for $x_0 > 0$, we have $V_0
(x_0) = \int_{-\infty}^{x_0} v_0 (x) dx = -\int_{x_0}^{-\infty}
v_0 (x) dx$ (because of (\ref{0masssys})), and a similar argument
shows that $|x_0 V_0(x)|$ approaches $0$ as $x_0 \rightarrow
+\infty$. Now assume that $(\alpha_1 , \alpha_2)$ is an interval
on which $V_0$ does not change sign (suppose, without loss of
generality, it is positive on this interval), and $V_0(\alpha_1)=
0$ and $V_0(\alpha_2) = 0$. Then $\alpha_1 V_0(\alpha_1)= 0$ and
$\alpha_2 V_0(\alpha_2) = 0$ (in the case $\alpha_1$ or $\alpha_2$
is $\pm\infty$, the aforementioned argument would work). In this
interval, then, we will have
$$\int_{\alpha_1}^{\alpha_2} |V_0 (x) dx| = \int_{\alpha_1}^{\alpha_2} V_0 (x) dx$$
which,by integration by parts, is equal to
$$ \alpha_2 V_0 (\alpha_2) - \alpha_1 V_0 (\alpha_1) - \int_{\alpha_1}^{\alpha_2} x v_0 (x) dx
\leq \int_{\alpha_1}^{\alpha_2} |x v_0 (x)| dx.$$ Now take
summation over all the intervals in the aforesaid form.
\end{proof}

Substituting $u$ with $v + \fe$ in (\ref{utassys}), we get:
\begin{equation}
\begin{split}
 &v_t +A v_x - B v_{xx}\\
  &= -(\fe_t +A \fe_x - B \fe_{xx} +
\Gamma(\fe, \fe)_x)\\
&+ \bold {O} (|v|^2 + |\fe||v| + |(\fe + v)(\fe + v)_x| + |\fe
+v|^3)_x\\
&=: \Psi(x,t) + \fcal(v,\fe)_x.
\end{split} \label{khatisys}
\end{equation}
Using (\ref{fiith}) and (\ref{coefficient}), we get
\begin{equation}
\begin{split}
\Psi(x,t) = &-(\fe_t +A \fe_x - B \fe_{xx} + \Gamma(\fe, \fe)_x) \\
= &-\sum_{i=1}^{n}\sum_{j\neq k}\Gamma^i_{jk} (\fe^j \fe^k)_x r_i
- \sum_{i \neq j} \Gamma^i_{jj}(\fe^j)^2_x r_i + \sum_{i \neq j}
b^i_j \fe^j_{xx} r_i.
\end{split} \label{termha}
\end{equation}
(The whole point is that, this way, we get rid of the terms
$(\fe^i)^2_x r_i$ and $\fe^i_{xx} r_i.$)

The Green function for the linear part of (\ref{khatisys}), i.e.,
for $v_t+Av_x-Bv_{xx},$ is
$$G(x, t; y)= \sum_{i=1}^n(4\pi t)^{-\frac 12}
e^{\frac{(x-y-a_it)^2}{4\beta_i t}}r_i l_i^t + R (x, t; y),$$ with
the remainder $R$,
  $$R (x, t; y) = \bold {O} ((1+t)^{-1}
\sum_{i=1}^n e^{\frac{(x-y-a_it)^2}{M t}})$$ (for proof see
\textbf{\cite{LZe}}).

Using Duhamel's principle,
\begin{equation}
\begin{aligned}
v(x, t) &= \int_{-\infty}^{+\infty} G(x,t;y) v_0(y) dy\\
&+ \int_{0}^{t} \int _{-\infty}^{+\infty} G(x,t-s;y) (\fcal(v,
\fe)_y(y,s)+ \Psi(y,s)) dy \, ds \\
 &= \int_{-\infty}^{+\infty} G_y(x,t;y) V_0(y) dy\\
&+ \int_{0}^{t} \int _{-\infty}^{+\infty} G_y(x,t-s;y) \fcal(v,
\fe)(y,s) dy \, ds \\
&+ \int_{0}^{t} \int _{-\infty}^{+\infty} G(x,t-s;y) \Psi(y,s) ds.
\end{aligned}
\end{equation}

\begin{theo}
Assume the above setting, $u = v + \fe$ and assume that
$|u_0|_{L^1},\, |u_0|_{L^{\infty}}, \, |x u_0(x)|_{L^1} \leq E_0$,
for sufficiently small $E_0$ (these inequalities translate into
similar ones for $v_0$ and $\fe_0$). Then,
\begin{equation} |v(\cdot, t)|_p \leq C \,E_0 (1+t)^{-\frac
12(1-\frac 1p) - \frac 14}
\end{equation}
  for some constant $C$.
\end{theo}
\begin{proof}
Let
\begin{equation}
\zeta(t) := \sup_{0\leq s \leq t, \, 1 \leq p \leq \infty
}|v(\cdot, s)|_{L^p}(1+s)^{\frac{1}{2}(1-\frac 1p)+\frac 14},
\end{equation}
 i.e.,
$|v(\cdot, s)|_{L^p} \leq (1+s)^{-\frac{1}{2}(1-\frac 1p)-\frac
14} \zeta(t)$, and in particular $|v(\cdot, s)|_{L^\infty} \leq
(1+s)^{-\frac 34} \zeta(t)$. The goal is to show:
 $$\zeta(t) \leq C(E_0 + \zeta(t)^2).$$ But
then, if $E_0$ is sufficiently small, this implies that $\zeta(t)
\leq 2C E_0$, and that is what we are looking for.

 Obviously
$|G_y|_{L^p} \leq \tpnim$. Also $|\fe(\cdot, t)|_{L^p} \leq E_0
(1+t)^{-\frac 12(1-\frac 1p)}$ and $|\fe_x(\cdot, t)|_{L^p} \leq
E_0 (1+t)^{-\frac 12(1-\frac 1p)-\frac12}$. We need to find some
bounds for $v_x$, and that is the subject of the following lemma:
\begin{lem}
Given above setting, we will have:
\begin{equation}
|v_x(\cdot,t)|_{L^p} \le \begin{cases}
C(|v(\cdot, t-1)|_{L^p} + |\mcal(\fe(\cdot, t)|_{L^p}) , &\hbox{for}\ t \ge 1\\
C (t^{-\frac{1}{2}}|v_0|_{L^p}|+|\mcal(\fe(\cdot, t)|_{L^p}) , & \hbox{for}\ t\le 1.\\
\end{cases} \label{vx1}
\end{equation}
where $\mcal(\fe):= -\fe_t + (B(\ubar +\fe)\fe_x)_x -
(f(\ubar+\fe))_x$ and consequently
\begin{equation}
|v_x(\cdot,t)|_{L^p}\leq \begin{cases}  C(E_0+\zeta(t))t^{-\frac
12(1-\frac1p)-\frac 14}&\hbox{for}\ t \ge 1\\
C (E_0+\zeta(t))t^{-\frac{1}{2}} , & \hbox{for}\ t\le 1.
\end{cases}
\label{vx2}
\end{equation}
\end{lem}
\begin{proof}
$$ \utild_t+f(\utild)_x=(B(\utild) \utild_x)_x$$
implies
$$v_t+ \fe_t +f(\ubar +\fe
+v)_x=(B(\ubar+\fe+v)(\ubar_x+\fe_x+v_x))_x.
$$ Therefore,
$$v_t+ (f(\ubar +\fe +v)-f(\ubar+\fe))_x-((B(\ubar +\fe +v)-B(\ubar +\fe))(\ubar_x+\fe_x))_x$$
$$-(B(\ubar +\fe +v)v_x)_x$$
$$=-\fe_t-f(\ubar+\fe)_x+(B(\ubar +\fe )\fe_x)_x$$
$$=:\mcal (\fe).$$
From here we can use short time estimates described in
\textbf{\cite{ZH}}, section 11, to achieve inequality (\ref{vx1})
(see the argument in lemma \ref{wxlem}).
 (\ref{vx2}) follows immediately, if we notice that, by
definition of $\fe$, we have $\mcal(\fe)=\bold O (|\fe||\fe_x|).$
\end{proof}
 Returning to the proof of the theorem, whenever $0\leq s \leq t,$ then,
\begin{equation}
|\fcal(v, \fe)(\cdot, s)|_{L^p} \leq C(E_0 +
\zeta^2(t))s^{-\frac12(1-\frac1p)-\frac34},\label{few5}
\end{equation}
(we used $E_0 \leq 1$ and so $E^2_0, E^3_0 \leq E_0$ and $E_0
\zeta (t) \leq \frac 12 E^2_0 + \frac 12 \zeta^2(t)$\,).

When $t \leq 1$, then
\begin{equation}
\begin{aligned} |v(\cdot,& t)|_{L^p}
\le C |v_0|_{L^p}\\
&+\int^t_0 |G_y|_{L^1}|\fcal(v, \fe)|_{L^p}(s)ds\\
&+\int^t_0 (|G|_{L^1}||\psi(y,s)|_{L^p}(s)ds\\
 &\le C E_0+(E_0+\zeta(t)^2) \int^t_0 (t-s)^{-\frac{1}{2}}
s^{-\frac{1}{2}}ds +C E_0 \int^t_0 (1+s)^{\frac 12}\\
&\le C (E_0 +\zeta(t)^2)\leq C (E_0 +\zeta(t)^2)(1+t)^{-\frac12(1-\frac1p)-\frac14}.\\
\end{aligned} \label{k1}
\end{equation}

 For $t \geq 1$ we use again Haussdorf-Young inequality
 to obtain:
 {\allowdisplaybreaks
\begin{equation}
\begin{aligned}
&\left|\int_{-\infty}^{+\infty} G_y(x,t;y) V_0(y) dy\right|_{L^P}\\
&\leq  |V_0|_{L^1} |G_y|_{L^p}\leq C E_0
t^{-\frac{1}{2}(1-1/p)-\frac 12},\\
\end{aligned} \label{b10}
\end{equation}
and
\begin{equation}
\begin{aligned}
|\int_{0}^{t}& \int _{-\infty}^{+\infty} G_y(x,t-s;y) \fcal(v,
\fe)(y,s) dy \, ds|_{L^p} \\
\leq &\int^{t/2}_0 |G_y|_{L^p}|\fcal(v, \fe)(\cdot, s)|_{L^1}ds\\
&+ \int^t_{t/2}|G_y|_{L^1}|\fcal(v, \fe)(\cdot, s)|_{L^p}ds\\
\leq &C (E_0 + \zeta (t)^2)(\int^{t/2}_0
(t-s)^{-\frac{1}{2}(1-1/p)-\frac{1}{2}}
s^{- \frac 34}ds\\
&+\int^t_{t/2}(t-s)^{-\frac{1}{2}}
s^{-\frac{1}{2}(1-1/p) - \frac 34} ds)\\
\leq &C(\zeta_0 +\zeta (t)^2)t^{-\frac{1}{2}(1-1/p) -\frac 14}.\\
\leq &2C(\zeta_0 +\zeta (t)^2)(1+t)^{-\frac{1}{2}(1-1/p) -\frac
14}.
\end{aligned} \label{b101}
\end{equation}
}
 It remains only  to deal with the term $\Psi$ in (\ref{khatisys}), which, by
(\ref{termha}), includes the terms in the form $\Gamma^i_{jk}
(\fe^j \fe^k)_x r_i$ for $j \neq k,$ and $\Gamma^i_{jj}
(\fe^j)^2_x r_i$ and $b^i_j \fe^j_{xx} r_i$ for $i \neq j.$

 Lemma \ref{interaction} takes care of the
terms in the form $(\fe^j \fe^k)$ for $j\neq k,$ as strict
hyperbolicity of $A$ implies $a_j \neq a_k.$, hence giving us
\begin{equation}
|(\fe^j \fe^k)|_{L_p} \leq C E^2_0 e^{-\eta s} \leq C E_0
(1+s)^{-\frac12(1-\frac1p)-\frac34}. \label{few6}
\end{equation}
This then will be treated similar to the way (\ref{few5}) is
treated in (\ref{b101}).

 For other terms, we need to estimate
$$\left| \int_{0}^{t} \int_{-\infty}^{+\infty} G(x, t; y)(\fe^j)^2_x
r_i \right|_{L^p}$$ for $i \neq j.$ As the remainder $R$ in $G$ is
small enough, the only part of concern would be:
$$\left| \int_{0}^{t} \int_{-\infty}^{+\infty} (4\pi t)^{-\frac 12}
e^{\frac{(x-y-a_it)^2}{4\beta_i t}} (\fe^j)^2_x  \right|_{L^p}$$
for $i\neq j$. Here we have a  heat kernel convecting at the speed
$a_i$ convoluted against  a diffusion wave $\fe^j$ which is
similar to $\phi$ in Example (\ref{phi}), but convecting at the
speed $a_j.$ Corollary (\ref{phicor}), then,  implies that the
above
term would be less than  $CE_0(1+t)^{-\frac 12(1-\frac1p)-\frac 14}$.\\
The term $b^i_j \fe^j_{xx} r_i$ can be treated similarly. This
finishes the proof.
\end{proof}

\begin{rem}
\textup{ In the scalar case, we get rid of the terms in the form
of $\fe^2$ all together, since $i=j=1$ is the only possibility.
Hence it can readily be seen that, in the scalar case, we would
obtain, using this argument, the decay rates: $|v(\cdot, t)|_{L^p}
\sim (1+t)^{-\frac 12(1-\frac 1p)-\frac 12 +\epsilon}$, for
$\epsilon$ arbitrarily small \textbf{\cite{Liu}}}.
\end{rem}
\section{Strictly parabolic cases with a viscous shock solution } \label{S:overc}
We now focus on a shock wave solution of the system of viscous
conservation laws
\begin{equation}
 \utild_t+f(\utild)_x=(B(\utild) \utild_x)_x,  \label{oveq}
\end{equation}
  where $f\in \mathbb{R}^n$ and $B(\utild) \in \mathbb{R}^{n\times n}$, and
  $\utild \in \mathbb{R}^n$ is a
 perturbation of (without loss of
generality) a stationary viscous shock solution
\begin{equation}
 \ubar=\ubar(x),  \, \lim_{x \to \pm \infty} \ubar(x)=:
u_\pm, \label{ovhad}
\end{equation}
i.e., $\ubar$ solves
\begin{equation}
\ubar'= B(\bar u)^{-1}(f(\bar u) - f(u_-)). \label{ovode}
\end{equation}
 Following \textbf{\cite{ZH, Z}}, we make assumptions
($\hcal$) below.

\medskip
\noindent Assumptions ($\hcal$):
\medskip

({$\hcal0$}) \quad  $f,B\in C^3$.
\medskip

({$\hcal1$}) \quad $Re \, \sigma(B) > 0$.
\medskip

({$\hcal2$}) \quad $\sigma (df(u_\pm))$ real, distinct, and
nonzero.
\medskip

({$\hcal3$}) \quad $Re\,  \sigma(-ik df(u_\pm) -k^2 B(u_\pm))<
-\theta k^2$ for all real $k$, some $\theta>0$.
\medskip

({$\hcal4$}) \quad All set of the stationary solutions near
$\ubar$ of (\ref{oveq})-(\ref{ovhad}), connecting the same values
$u_\pm$ forms a smooth  manifold $\{\ubar^{\delta}\}, \delta\in
\mathbb{R}^{\ell}, \ubar^0=\ubar.$ Moreover the stable manifold of
$u_-$ and the unstable manifold of $u_+$ (with respect to
(\ref{ovode})) are transverse.

 Condition $(\hcal3)$ is the {\it stable viscosity
matrix} criterion of Majda and Pego, corresponding to linearized
stability of the constant solutions $u \equiv u_\pm$
\textbf{\cite{MP,Kaw}} (clearly necessary for stability of
$\ubar(\cdot)$ of the type we seek, see further discussion
(\textbf{\cite{ZH}}, pp. 746, 767, and 774--775). Note that
condition $(\hcal4)$ is the condition $H4$ of \textbf{\cite{ZH}},
plus the assertion that the shock is of ``standard" or ``pure"
type (see \textbf{\cite{ZH}}, section 10). This implies that we
have $n+\ell$ incoming characteristics, entering the shock, hence
$n-\ell$ outgoing modes, i.e., eigenvalues are in the form:
\begin{equation}
a_1^- < \cdots <a_{p-1}^-< 0 < a_p^- <\cdots <a_n^-, \label{pchar}
\end{equation}
 and
\begin{equation}
a_1^+ < \cdots <a_{p+\ell -1}^+< 0 < a_{p+\ell}^+ <\cdots <a_n^+,
\label{order}
\end{equation}
where $a_i^\pm$ denote the (ordered) eigenvalues of $df(u_\pm)$.
If $\ell = 1$, we have a Lax type shock wave, in which case there
are $n-1$ outgoing modes (corresponding to $a_i^{\pm}\gtrless 0$)
and $n+1$ incoming modes (corresponding to $a_i^{\pm}\lessgtr 0$).
If $\ell > 1$, we have an overcompressive shock, with $n-\ell$
outgoing modes and $n+\ell$ incoming modes. For further discussion
see \textbf{\cite{ZH}}.

The following Lemma proved in \textbf{\cite{MP}}  asserts that
$u_\pm$ are hyperbolic also in the ODE sense  (for an alternative
proof, see Remark 2.3 in section 2). This implies exponential
approach of $\bar u^\delta$ to its asymptotic states at
$x=\pm\infty$, a fact that will be crucial in our subsequent
analysis. See \textbf{\cite{MP}} and also \textbf{\cite{ZH}} for
proofs.
 \begin{lem}
\label{mplem} Given $(\hcal0)-(\hcal3)$, the stable/unstable
manifolds of $df(u_\pm)$ and $B(u_\pm) ^{-1}df(u_\pm)$ have equal
dimensions. In particular,  $B(u_\pm) ^{-1}df(u_\pm)$ has no
center manifold.
\end{lem}
\begin{cor} \label{ubardecay} Given  $(\hcal0)-(\hcal4)$,
solutions $\bar u^\delta$ of (\ref{ovode}) are in $C^{4}$,
satisfying
$$D_x^j D_\delta^i(\bar u^\delta(x)-u_{\pm})=\bold{O}(e^{-\alpha|x|}),
\quad \alpha>0, \, 0\le j\le 4, \,i=0,1, $$ as $x\rightarrow
\pm\infty$
\end{cor}

Linearizing about $\ubar(\cdot)$ gives:
\begin{equation}
v_t=Lv:=-(Av)_x+(Bv_x)_x, \label{linearov}
\end{equation}
with
\begin{equation}
B(x):= B(\ubar(x)), \quad  A(x)v:=
df(\ubar(x))v-dB(\ubar(x))v\ubar_x. \label{AandBov}
\end{equation}
Denoting $A^\pm := A(\pm \infty)$,  $B^\pm:= B(\pm \infty)$, and
considering  corollary \ref{ubardecay}, it follows that
\begin{equation}
|A(x)-A^-|= \bold {O} (e^{-\eta |x|}), \quad |B(x)-B^-|= \bold {O}
(e^{-\eta |x|}) \label{ABboundsov}
\end{equation}
as $x\to -\infty,$ for some positive $\eta.$ Similarly for $A^+$
and $B^+,$ as $x\to +\infty.$ Also $|A(x)-A^\pm|$ and
$|B(x)-B^\pm|$ are bounded for all $x$.

 Define the {\it (scalar) characteristic speeds} $a^\pm_1<
\cdots < a_n^\pm$ (as above) to be the eigenvalues of $A^\pm$, and
the left and right {\it (scalar) characteristic modes} $l_j^\pm$,
$r_j^\pm$ to be corresponding left and right eigenvectors,
respectively (i.e., $A^\pm r_j^\pm = a_j^\pm r_j^\pm,$ etc.),
normalized so that $l^+_j \cdot r^+_k=\delta^j_k$ and $l^-_j \cdot
r^-_k=\delta^j_k$. Following Kawashima \textbf{\cite{Kaw}}, define
associated {\it effective scalar diffusion rates}
$\beta^\pm_j:j=1,\cdots,n$ by relation
\begin{equation}
\left(
\begin{matrix}
\beta_1^\pm &&0\\
&\vdots &\\
0&&\beta_n^\pm
\end{matrix}
\right) \quad = \hbox{diag}\ L^\pm B^\pm R^\pm, \label{betaov}
\end{equation}
where $L^\pm:=(l_1^\pm,\dots,l_n^\pm)^t$,
$R^\pm:=(r_1^\pm,\dots,r_n^\pm)$ diagonalize $A^\pm$.

Let
\begin{equation}
G(x,t;y):= e^{Lt}\delta_y (x) \label{ov3.7}
\end{equation}
 be the Green's function associated with $(\partial_t -L)$. Then,
the relevant linearized theory can be summarized in the following
two propositions, proved in \textbf{\cite{ZH}}.

\begin{prop} \label{Dov}
  Given ($\hcal$), necessary conditions for
$L^p$-linearized orbital stability, $p>0$, of $\bar u(\cdot)$ with
respect to perturbations $v_0\in C^\infty_0$ are:

\noindent Assumptions ($\dcal$):

$(\dcal1)$ \quad  $L$ has no ($L^2$, without loss of generality)
eigenvalues in $\{Re \lambda \ge 0\} \setminus \{0\}$.
\medskip

$(\dcal2)$ \quad $\{r^\pm; a^\pm \gtrless 0\} \cup
\{\int_{-\infty}^{+\infty}\frac{\partial \ubar^\delta}{\partial
\delta_i} dx; i=1, \cdots, \ell \}$ is a basis for $\mathbb{R}^n$,
with $\int_{-\infty}^{+\infty}\frac{\partial
\ubar^\delta}{\partial \delta_i} dx$  computed at $\delta=0.$
\end{prop}
\begin{prop} \label{greenovc}
  Under assumptions  ($\hcal$), ($\dcal$), we have for $y\le
0$ the decomposition
\begin{equation} G=E+S+R,\label{GESRov}\end{equation}
where
\begin{equation}
\begin{aligned} E(x,t;y) &:= \sum_{i=1}^{
\ell}\sum_{a_k^- > 0} [c^{0,i}_{k,-}]\frac{\partial\ubar^
\delta}{\partial \delta_i} (x) {l_k^-}^{t}
\left(errfn \left(\frac{y+a_k^-t}{\sqrt{4\beta_k^-t}}\right)\right.\\
&-errfn\left.\left(\frac{y-a_k^-t}{\sqrt{4\beta_k^-t}}\right)\right),\\
\end{aligned}
\label{Eov}
\end{equation}
\begin{equation}
\begin{aligned} S(x,t;y)&:= \chi_{\{t\ge 1\}}\sum_{a_k^-<0}r_k^-  {l_k^-}^t
(4\pi \beta_k^-t)^{-1/2} e^{-(x-y-a_k^-t)^2 / 4\beta_k^-t} \\
&+\chi_{\{t\ge 1\}} \sum_{a_k^- > 0} r_k^-  {l_k^-}^t (4\pi
\beta_k^-t)^{-1/2} e^{-(x-y-a_k^-t)^2 / 4\beta_k^-t}
\left({\frac {e^{-x}}{e^x+e^{-x}}}\right)\\
&+ \chi_{\{t\ge 1\}}\sum_{a_k^- > 0, \,  a_j^- < 0}
[c^{j,-}_{k,-}]r_j^-  {l_k^-}^t (4\pi \bar\beta_{jk}^- t)^{-1/2}
e^{-(x-z_{jk}^-)^2 / 4\bar\beta_{jk}^- t}
\left({\frac{e^{ -x}}{e^x+e^{-x}}}\right),\\
&+\chi_{\{t\ge 1\}} \sum_{a_k^- > 0, \,  a_j^+ > 0}
[c^{j,+}_{k,-}]r_j^+  {l_k^-}^t (4\pi \bar\beta_{jk}^+ t)^{-1/2}
e^{-(x-z_{jk}^+)^2 / 4\bar\beta_{jk}^+ t}
\left({\frac{e^{ x}}{e^x+e^{-x}}}\right),\\
\end{aligned}
\label{Sov}
\end{equation}
with
\begin{equation}
z_{jk}^\pm(y,t):=a_j^\pm\left(t-\frac{|y|}{|a_k^-|}\right)
\label{zjkov}
\end{equation}
and
\begin{equation}
\bar \beta^\pm_{jk}(x,t;y):= \frac{x^\pm}{a_j^\pm t} \beta_j^\pm +
\frac{|y|}{|a_k^- t|} \left( \frac{a_j^\pm}{a_k^-}\right)^2
\beta_k^-, \label{betaaverageov}
\end{equation}
and
\begin{equation}
\begin{aligned} R&(x,t;y)=\\
&\bold{O}(e^{-\eta(|x-y|+t)})\\
&+\sum_{k} \bold {O} \left( (t+1)^{-1/2}
e^{-\eta x^+} +e^{-\eta|x|} \right)
t^{-1/2}e^{-(x-y-a_k^- t)^2/Mt} \\
&+ \sum_{a_k^- > 0, \, a_j^- < 0} \chi_{\{ |a_k^- t|\ge |y| \}}
\bold {O} ((t+1)^{-1/2}t^{-1/2} e^{-(x-a_j^-(t-|y/a_k^-|))^2/Mt}
e^{-\eta x^+}, \\
&+ \sum_{a_k^- > 0, \, a_j^+> 0} \chi_{\{ |a_k^- t|\ge |y| \}}
\bold {O} ((t+1)^{-1/2}t^{-1/2} e^{-(x-a_j^+ (t-|y/a_k^-|))^2/Mt}
e^{-\eta x^-}, \\
\end{aligned}
\label{Rov}
\end{equation}
for some $\eta$, $M>0$, where $x^\pm$ denotes the
positive/negative part of $x$, indicator function $\chi_{\{ |a_k^-
t|\ge |y| \}}$ is one for $|a_k^- t|\ge |y|$ and zero otherwise,
indicator function $\chi_{\{ t\ge 1 \}}$ is one for $t\ge 1$ and
zero otherwise, and {\it scattering coefficients}
$[c_{k,-}^{0,i}]$, $[c_{k,-}^{j,\pm}]$ are constant, with
\begin{equation}
\sum_{a_j^- < 0} [c_{k, \, -}^{j, \, -}]r_j^- + \sum_{a_j^+ > 0}
[c_{k, \, -}^{j, \, +}]r_j^+ + \sum_{i=1}^{\ell}[c_{k,-}^{0,
i}]\int_{-\infty}^{+\infty}\frac{\partial \ubar^\delta}{\partial
\delta_i} dx = r_k^- \label{scatteringov}
\end{equation}
for each $k$ (note: uniquely determined, by condition (D2)), and
\begin{equation}
\begin{aligned} \sum_{a_k^->0} [c_{k,-}^{0,i}] l_k^- &=
 \sum_{a_k^+<0} [c_{k,+}^{0,i}] l_k^+ \\
&=\pi_i := (r_1^-,\dots,r_{p-1}^-,r_{p+l}^+,\dots,r_n^+,
\int_{-\infty}^{+\infty}\frac{\partial \ubar^\delta}{\partial
\delta_i} dx)^{-1}e_{n-i+1},
\end{aligned}
\label{pi}
\end{equation}
where $e_j$ denotes the $j$th standard basis element, and with
$\frac{\partial \ubar^\delta}{\partial \delta_i} $ always computed
at $\delta=0.$ Likewise, we have the derivative bounds
\begin{equation}
\begin{aligned} |&R_x|=\\
&\bold{O}(e^{-\eta(|x-y|+t)})\\
&+\sum_{k} \bold {O} \left( (t+1)^{-1/2}t^{-1/2}
e^{-\eta x^+} +e^{-\eta|x|} \right)
t^{-1/2}e^{-(x-y-a_k^- t)^2/Mt} \\
&+ \sum_{a_k^- > 0, \, a_j^- < 0} \chi_{\{ |a_k^- t|\ge |y| \}}
\bold {O} ((t+1)^{-1}t^{-1/2} e^{-(x-a_j^-(t-|y/a_k^-|))^2/Mt}
e^{-\eta x^+}, \\
&+ \sum_{a_k^- > 0, \, a_j^+ > 0} \chi_{\{ |a_k^- t|\ge |y| \}}
\bold {O} ((t+1)^{-1}t^{-1/2} e^{-(x-a_j^+(t-|y/a_k^-|))^2/Mt}
e^{-\eta x^-}, \\
\end{aligned}
\label{Rx}
\end{equation}
\begin{equation}
\begin{aligned} |&R_y|=\\
&\bold{O}(e^{-\eta(|x-y|+t)})\\
&+ \sum_{k} \bold {O} \left( (t+1)^{-1/2} e^{-\eta
x^+} +e^{-\eta|x|} \right)
t^{-1}e^{-(x-y-a_k^- t)^2/Mt} \\
&+ \sum_{a_k^- > 0, \, a_j^- < 0} \chi_{\{ |a_k^- t|\ge |y| \}}
\bold {O} ((t+1)^{-1}t^{-1/2} e^{-(x-a_j^-(t-|y/a_k^-|))^2/Mt}
e^{-\eta x^+}. \\
&+ \sum_{a_k^- > 0, \, a_j^+ > 0} \chi_{\{ |a_k^- t|\ge |y| \}}
\bold {O} ((t+1)^{-1}t^{-1/2} e^{-(x-a_j^+(t-|y/a_k^-|))^2/Mt}
e^{-\eta x^-}. \\
\end{aligned}
\label{Ry}
\end{equation}
A symmetric decomposition holds for $y\ge0$.
Moreover, for $|x-y|/t$ sufficiently large,
\begin{equation}\label{gaussian}
|G|\le C e^{-\frac{|x-y|^2}{Mt}}.
\end{equation}
\end{prop}
\begin{rem}
\textup{ Though it was not remarked in {\bf \cite{MaZ.3}}, the
terms $E$ and $S$ are continuous at $y=0$, a consequence of the
respective scattering relations \eqref{pi} and
\eqref{scatteringov}. (Note that values at $y=0$ correspond to
time-asymptotic states described by the scattering relations,
which depend only on mass and not position of data.) }
\end{rem}

\begin{rem}\textup{
The term $e^{-\eta(|x-y|+t)}$ in $R$ and its derivatives corrects
a minor omission in \textbf{\cite{Z}}. This term comes from the
fact that, in the far field, $E$ and $S$ decay at this rate while
entire $G$ decays at faster Gaussian rate. The Gaussian decay
\eqref{gaussian} was proved but not stated in
\textbf{\cite{MaZ.3}}. The bound for $R_x$ is given here only for
the sake of completeness, and is not going to play a role in our
calculations.}
\end{rem}

\begin{rem}
\textup{In \textbf{\cite{Z}} and \textbf{\cite{MaZ.3}} the above
bounds have been explicitly stated and proved for Lax case, and
only some hints are given as about the overcompressive case. The
proof for the overcompressive case, however, is not very different
and can be achieved following the same outline given for Lax
case.}
\end{rem}

 Define $e_i, i=1, \cdots, \ell$ for $y < 0$
\begin{equation}
\begin{aligned} e_i(y,t) &:= \sum_{a_k^- > 0} [c^{0,i}_{k,-}]
{l_k^-}^t
\left(errfn \left(\frac{y+a_k^-t}{\sqrt{4\beta_k^-t}}\right)\right.\\
&-errfn\left.\left(\frac{y-a_k^-t}{\sqrt{4\beta_k^-t}}\right)\right),\\
\end{aligned}
\label{eov}
\end{equation}
Hence
\begin{equation}
E(x,t:y)=\sum_{i=1}^{\ell}\frac{\partial\ubar^\delta}{\partial
\delta_i} (x)e_i(y,t)
\end{equation}
 and symmetrically for $y > 0.$ Define also
\begin{equation}
\gtild = S+R.
\end{equation}

We have the following bounds for $\tilde G$ and $e_i$'s:
\begin{lem} \label{gebounds}
  Under assumptions $(\hcal)$ and $(\dcal)$ there holds
\begin{equation}
|\int_{-\infty}^{+\infty} \tilde G(\cdot,t;y)f(y)dy|_{L^p} \le C
\min \{|f|_{L^p}, t^{-\frac{1}{2}(1-1/p)} |f|_{L^1} \},
\label{Gov}
\end{equation}
\begin{equation}
|\int_{-\infty}^{+\infty} \tilde G_y(\cdot,t;y)f(y)dy|_{L^p} \le C
\min \{t^{-1/2}|f|_{L^p}, t^{-\frac{1}{2}(1-1/p)-1/2} |f|_{L^1}
\}, \label{Gyov}
\end{equation}
for all $t\ge 0$, $f\in L^1\cap L^p$, some $C>0$.
\end{lem}
\begin{lem} \label{elemov}
 The kernels $e_i$'s satisfy
\begin{equation}
|e_{i_y} (\cdot, t)|_{L^p},  |e_{i_t}(\cdot, t)|_{L^p} \le C
t^{-\frac{1}{2}(1-1/p)}, \label{eypov}
\end{equation}
\begin{equation}
|e_{i_{ty}}(\cdot, t)|_{L^p} \le C t^{-\frac{1}{2}(1-1/p)-1/2},
\label{eytpov}
\end{equation}
for all $t>0$.  Moreover, for $y\leq 0$ we have the pointwise
bounds
$$
|e_{i_y} (y,t)|, |e_{i_t} (y,t)| \le C \sum_{a^-_k
> 0}t^{-\frac{1}{2}} e^{-\frac{(y+a^-_kt)^2}{Mt}},
$$
$$
|e_{i_{ty}} (y,t)| \le C \sum_{a^-_k
> 0}t^{-1} e^{-\frac{(y+a^-_kt)^2}{Mt}},
$$
for $M>0$ sufficiently large (i.e., $>4b_\pm$), and symmetrically
for $y\ge 0$.
\end{lem}
\begin{proof}
See \textbf{\cite{Z}} and \textbf{\cite{MaZ.4}}.
\end{proof}

Let $\utild$ solve (\ref{oveq}), and, using (D2), assume that
$$\int_{-\infty}^{+\infty} \utild (x, 0) - \ubar(x)= \sum_{a_j^- <0}m_j r_j^- +
\sum_{a_j^+ >0}m_j r_j^+ +\sum_{i=1}^{\ell}\int c_i\frac{\partial
\ubar^\delta}{\partial \delta_i}
$$
with $m_i$'s and $c_i$'s small enough. Using the Implicit Function
Theorem, we can find $\delta_0$ such that
$$\int_{-\infty}^{+\infty} \utild (x, 0) - \ubar^{\delta_0}(x)= \sum_{a_j^- <0}m'_j r_j^- +
\sum_{a_j^+ >0}m'_j r_j^+
$$
where each $m'_i$ is just ``slightly" different from $m_i$. Notice
that this way we have no ``mass" in any $\int \frac{\partial
\ubar^\delta}{\partial \delta_i}$ direction anymore. Therefore, by
replacing $\ubar$ with $\ubar^{\delta_0}$ and without loss of
generality, we can assume $\delta_0=0$ and
$$\int_{-\infty}^{+\infty} \utild (x, 0) - \ubar(x)= \sum_{a_j^- <0}m_j r_j^- +
\sum_{a_j^+ >0}m_j r_j^+.$$
\begin{rem}\textup{
In Lax case shock waves, $\ubar^\delta(x)=\ubar(x+\delta),$ hence
$\ddp=u'(x),$ and $\delta_0$ can be explicitly computed:
$\delta_0=c_1.$ }
\end{rem}
Let $u(x,t)=\utild(x,t)-\ubar(x)$ and use Taylor's expansion
around $\ubar^{\delta(t)}(x)$ to find
\begin{equation}
u_t + (A(x)u)_x - (B(x)u_x)_x = -(\Gamma(x)(u,u))_x + Q(u, u_x)_x,
\label{utaylorov}
\end{equation}
where $\Gamma(x)(u,u) = d^2f(\ubar)(u,u)- d^2B(\ubar)(u,u)\ubar_x$
and $$Q(u, u_x) = \bold {O} (|u||u_x|+|u|^3).$$  Denote
$\Gamma^\pm = \Gamma(\pm \infty),$ and note that we have similar
statements to (\ref{ABboundsov}) for $\Gamma(x)-\Gamma^\pm.$
Define constant coefficients $b^{\pm}_{ij}$ and
$\Gamma_{ijk}^{\pm}$ to satisfy
\begin{equation}
\Gamma^\pm (r^\pm_j, r^\pm_k) = \sum_{i=1}^n \Gamma_{ijk}^{\pm}
r^\pm_i,\quad B^\pm r^\pm_j= \sum_{i=1}^n b^{\pm}_{ij} r^\pm_i
\label{coefshock}
\end{equation}
hence of course $\beta^\pm_i = b^{\pm}_{ii}$, and  denote
$\gamma^\pm_i := \Gamma_{iii}^{\pm}.$
\begin{rem}
\textup{As it is pointed out in \textbf{\cite{Liu, Liu85, Liu97}},
 $\gamma^\pm_i \ne 0 $ is equivalent to \emph{genuine nonlinearity} of
the $i^{\textup{th}}$ field, and  $\gamma^\pm_i = 0 $ means that
the $i^{\textup{th}}$ field is \emph{linearly degenerate}.}
\end{rem}

We define diffusion waves along outgoing modes: for $a_i^- <0$
define the diffusion wave $\fe_i$ to solve:
\begin{equation}
\begin{cases}
\fe^i_t + a^-_i\fe^i_x - \beta^-_i \fe^i_{xx} = -\gamma^-_i
(\fe^{i^2})_x \qquad
\text{for} &t > -1\\
\fe^i (x, -1) = m_i\delta_0 &t= -1
\end{cases} \label{ovfiithshock-}
\end{equation}
and likewise for $a_i^+>0 $, define $\fe_i$ to be the solution of
\begin{equation}
\begin{cases}
\fe^i_t + a^+_i\fe^i_x - \beta^+_i \fe^i_{xx} = -\gamma^+_i
(\fe^{i^2})_x \qquad
\text{for} &t > -1\\
\fe^i (x, -1) = m_i\delta_0 &t= -1
\end{cases} \label{ovfiithshock+}
\end{equation}
and set
$$\fe =\sum_{a_i^- <0}\fe^i r^-_i + \sum_{a_i^+>0}^{n}\fe^i r^+_i .$$
Let $v:=u-\fe- \frac{\partial\ubar^\delta}{\partial
\delta}\delta(t),$ where
$\delta(t)=(\delta_1(t),\cdots,\delta_\ell(t))^{\textup{tr}}$ is
to be defined later, assuming $\delta(0) = 0.$ Notice that
\begin{equation}
\int_{-\infty}^{+\infty} v(x, 0) dx = 0, \label{zmassov}
\end{equation}
so if $V_0=\int_{-\infty}^{x} v(y,0) dy$ then by lemma \ref{wW}
$V_0 \in L^1$ and $|V_0|_{L^1} \leq |x v_0|_{L^1}$.

 Replacing $u$
with $v + \fe+\ddp\delta(t)$ in (\ref{utaylorov}) (
$\frac{\partial \ubar^\delta}{\partial \delta_i} $  computed at
$\delta=\delta_0=0$), and using the fact that
$\frac{\partial\ubar^\delta}{\partial \delta_i}$ satisfies the
linear time independent equation $Lv=0$, we will have
\begin{equation}
v_t - Lv = \Psi(x,t)+ \fcal(\fe, v,  \ddp\delta(t))_x + \ddp\dot
\delta (t), \label{khatiov}
\end{equation}
where
\begin{equation}
\begin{split}
\fcal(\fe, v, \ddp\delta &) = \bold {O}  ( |v|^2 + |\fe||v|
+ |v||\ddp\delta|+ |\fe|| \ddp\delta|+|\ddp\delta|^2 \\
&+ |(\fe + v+\ddp\delta)(\fe + v + \ddp\delta)_x| + |\fe
+v+\ddp\delta|^3 ).
\end{split}
\end{equation}
and $\Psi := -\fe_t - (A(x)\fe)_x  + (B(x)\fe_x)_x -(\Gamma(x)(\fe
, \fe))_x$. For $\Psi$  we write
\begin{equation}
\begin{split}
\Psi(x,t) = &-(\fe_t +A \fe_x - B \fe_{xx} + \Gamma(\fe, \fe)_x) \\
= &-\sum_{a_i^- < 0} \fe_t^i r_i^- + (A(x)\fe^i r_i^-)_x -
(B(x)\fe_x^i r_i^-)_x + (\Gamma(x)(\fe^ir_i^-,\fe^ir_i^-))_x\\
&- \sum_{a_i^+ > 0} \fe_t^i r_i^+ + (A(x)\fe^i r_i^+)_x -
(B(x)\fe_x^i r_i^+)_x + (\Gamma(x)(\fe^ir_i^+,\fe^ir_i^+))_x\\
& - \sum_{i\neq j}(\fe_i \fe_j \Gamma(x)(r_i^\pm, r_j^\pm))_x.
\end{split} \label{ghatiov}
\end{equation}
Let us write a typical term of the first summation ($a_i^- < 0$)in
the following form:
\begin{equation}
\begin{split}
\fe_t^i & r_i^- + (A(x)\fe^i r_i^-)_x -
(B(x)\fe_x^i r_i^-)_x + (\Gamma(x)(\fe^ir_i^-,\fe^ir_i^-))_x\\
&= \left[(A(x)-A^-)\fe^i r_i^- -
(B(x)-B^-)\fe_x^i r_i^- + (\Gamma(x)-\Gamma^-)(\fe^ir_i^-,\fe^ir_i^-)\right]_x\\
&+ \fe_t^i r_i^- + (\fe^i_xA^- r_i^-) - (\fe_{xx}^i B^- r_i^-) +
((\fe^i)_x\Gamma^-(r_i^-,r_i^-)).
\end{split} \label{ghatitarov}
\end{equation}
Now we use the definition of $\fe^i$ in (\ref{ovfiithshock-}) and
the definition of coefficients $b_{ij}$ and $\Gamma_{ijk}$ in
(\ref{coefshock}) to write the last part of (\ref{ghatitarov}) in
the following form:
\begin{equation}
\begin{split}
\fe_t^i r_i^- + (\fe^i_xA^- r_i^-) - &(\fe_{xx}^i B^- r_i^-) +
((\fe^i)_x\Gamma^-(r_i^-,r_i^-))\\
&= -\fe^i_{xx}\sum_{j\ne i} b^-_{ij}r^-_j -
(\fe^i)^2_{x}\sum_{j\ne i} \Gamma^-_{jii}r^-_j.
\end{split} \label{ajabaaov}
\end{equation}
Similar statements hold for $a_i^+ > 0$ with minus signs replaced
with plus signs.

Later we will need some estimates for $v_x$. Short time estimates
gives us the necessary bounds. Let us first provide the requisite
short time existence/regularity theory for general quasilinear
parabolic systems, using the paramatrix method of Levi
\textbf{\cite{LSU, Le}}.
\begin{prop} \label{levi}
Let $\hat{A}(x,t)$, $\hat B(x,t)$, and $\hat C(x,t)$ be uniformly
bounded in $L^\infty$ and  $C^{(0,0)+(\alpha,\alpha/2)}(x,t)$,
$0<\alpha, <1$, taking values on a compact set, with $Re \, \sigma
(\hat B)$ positive and bounded strictly away from zero. Then, for
$0<t<T$, $T$ sufficiently small, there is a Green's function $\hat
G(x,t;y,s)\in C^{1,0}(x,t)$ associated with the Cauchy problem for
\begin{equation}
 v_t = \hat Cv +  (\hat Av)_x + (\hat Bv_x)_x, \quad  v\in \Bbb{R}^n,
\label{divergence}
\end{equation}
satisfying bounds
\begin{equation}
|D^j_x \hat G(x,t;y,s)| \le Ct^{-(j+1)/2}e^{-(x-y)^2/M(t-s)},
\quad j=0,1, \label{parametrix}
\end{equation}
where $C$, $M$, $T>0$ depend only on the bounds on the
coefficients and on the lower bound on $Re \, \sigma(\hat B)$.
\end{prop}
\begin{proof}
 See \textbf{\cite{ZH, LSU}}.
\end{proof}

The following lemma provides  us  with the short time estimates we
need for $v_x.$

\begin{lem} \label{wxlem}
 Given the above setting, and assuming $v$ remains bounded for all the time, we will have:
\begin{equation}
    |v_x(\cdot,t)|_{L^p} \le
 \begin{cases}
   C(|v(\cdot, t-1)|_{L^p} + |\mcal(\fe,\ddp\delta)|_{L^p}) , &\hbox{for}\ t \ge 1\\
    C (t^{-\frac{1}{2}}|v_0|_{L^p}|+|\mcal(\fe, \ddp\delta)|_{L^p}) , & \hbox{for}\ t\le 1.\\
 \end{cases} \label{wx1s}
\end{equation}
where
\begin{equation}
\mcal (\fe, \ddp\delta)=-\ddp\dot\delta + \bold O (|\fe
+\ddp\delta||\fe_x + (\ddp\delta)_x|+ e^{-k|x|}|\fe|).
\label{mcal}
\end{equation}
for some $k>0$
\end{lem}
\begin{proof}
$ \utild_t+f(\utild)_x=(B(\utild) \utild_x)_x$ and
$f(\ubar)_x=(B(\ubar)\ubar_x)_x$ implies
\begin{equation}
\begin{split}
v_t+ \fe_t &+\ddp\dot\delta +f(\ubar +\fe+\ddp\delta
+v)_x-f(\ubar)_x\\
&=(B(\ubar+\fe+\ddp\delta+v)(\ubar_x+\fe_x+(\ddp)_x\delta +v_x))_x
- (B(\ubar)\ubar_x)_x.
\end{split}
\end{equation}
Hence:
\begin{equation}
\begin{split}
& v_t+ (f(\ubar +\fe+\ddp\delta
+v)-f(\ubar+\fe+\ddp\delta))_x\\
-\Big((B(\ubar +&\fe+\ddp\delta +v)-B(\ubar +\fe+\ddp\delta))
(\ubar_x+\fe_x+(\ddp\delta)_x)\Big)_x\\
-(B(\ubar +&\fe+\ddp\delta +v)v_x)_x\\
&=-\fe_t-\ddp\dot\delta-(f(\ubar+\fe+\ddp\delta)-f(\ubar))_x\\
+((B(\ubar +&\fe+\ddp\delta )-B(\ubar ))\ubar_x)_x\\
+(B(\ubar
+&\fe+\ddp\delta )(\fe_x+(\ddp)_x\delta)_x\\
&=:\mcal (\fe, \ddp\delta).
\end{split} \label{ABM}
\end{equation}
It is not difficult to observe that
\begin{equation}
\begin{aligned}
\mcal (\fe, \ddp\delta)= -\fe_t -(A(x)(\fe +\ddp\delta))_x + (B(x)
(\fe_x+\frac{\partial\ubar^\delta_x}{\partial\delta}\delta)_x
-\ddp\dot\delta\\ + \bold O (|\fe +\ddp\delta||\fe_x +
(\ddp\delta)_x|.
\end{aligned}\label{mcal1}
\end{equation}

 Using (\ref{ghatiov}),
(\ref{ghatitarov}), (\ref{ABboundsov}) and the fact that
$L\ddp=0$, we conclude (\ref{mcal}). We use $$f(\eta +v) -f(\eta)=
\int_0^1df(\eta+\theta v)d\theta v$$ and similar equation for $B$
in order to write (\ref{ABM}) in the form
\begin{equation}
v_t +(\hat A(x,t) v)_x -(\hat B(x,t)v_x)_x =\mcal (\fe,
\ddp\delta).\label{hatab}
\end{equation}
Now $\hat A$ and $\hat B$ depend on $v$. Momentarily assume $v$
(hence $\hat A$ and $\hat B$) is in
$C^{(0,0)+(\alpha,\alpha/2)}(x,t)$.
\begin{equation}
\begin{aligned} v_x(x,t) &= \int_{-\infty}^{+\infty} \hat
G_x(x,t;y,t-T)v(y,t-T)\, dy \cr &\quad + \int_{t-T}^{t}
\int_{-\infty}^{+\infty} \hat G_x(x,t;y,s) \mcal(y,s) \, dy \, ds.
\end{aligned} \label{duhamel}
\end{equation}
By Duhamel's principle, where  $\hat G$ is the  Green's function
for (\ref{hatab}), and using the $\hat G_x$ bounds of proposition
\ref{levi} for divergence-form operators, we find that
$$
|v_x(\cdot, t)|_{\infty}\le C(T)^{-1/2}|v(\cdot, t-T)|_{\infty} +
C(T)^{1/2} |\mcal|_{\infty}.
$$
In particular, for $t\ge T$, we obtain a uniform H\"older (indeed,
Lipshitz) bound on $v(\cdot, t)$ depending only on the $L^\infty$
norm of $v(\cdot, t-T)$. By the (standard) method of extension, we
thus obtain uniform  H\"older continuity of $v$ so long as $|v|$
remains bounded.
 (\ref{wx1s})
follows using (\ref{duhamel}) and taking (without loss of
generality) $T=1$. For a more detailed discussion see
\textbf{\cite{ZH}} (section 11).
\end{proof}

Now we employ Duhamel's principle to get from  (\ref{khatiov}) :
\begin{equation}
\begin{aligned} v(x,t)
&=\int^{+\infty}_{-\infty} G(x,t;y)v_0(y)dy\\
&+\int^t_0
\int^{+\infty}_{-\infty}G(x,t-s;y)\fcal(\fe, v, \ddp\delta(s))_y(y,s) dy \, ds,\\
&+ \int^t_0\int^{+\infty}_{-\infty}G(x,t-s;y)\Psi(y,s) dy
\, ds \\
&+  \delta(t)\cdot\ddp\\
&=-\int^{+\infty}_{-\infty} G_y(x,t;y) V_0(y)dy\\
&-\int^t_0
\int^{+\infty}_{-\infty}G_y(x,t-s;y)\fcal(\fe, v, \ddp\delta(s))(y,s) dy \, ds,\\
&+ \int^t_0\int^{+\infty}_{-\infty}G(x,t-s;y)\Psi(y,s) dy
\, ds \\
&+\ddp\delta(t).\\
 \label{shswgreenov}
\end{aligned}
\end{equation}
(the last part of the above equation follows from
$\int_{-\infty}^{+\infty}
G(x,t;y)\frac{\partial\ubar^\delta}{\partial
\delta_i}(y)dy=e^{Lt}\frac{\partial\ubar^\delta}{\partial
\delta_i}=\frac{\partial\ubar^\delta}{\partial \delta_i}$ and
$\delta(0)=0$). Set
\begin{equation}
\begin{aligned}
 \delta_i (t)
&=\int^\infty_{-\infty}e_{i_y}(y,t) V_0(y)dy\\
&+\int^t_0\int^{+\infty}_{-\infty} e_{i_y}(y,t-s)\fcal(\fe,
v,\ddp\delta(s))
 dy ds\\
&- \int^t_0\int^{+\infty}_{-\infty}e_i(x,t-s;y)\Psi(y,s)dy \, ds.
\end{aligned} \label{deltaov}
\end{equation}
Using (\ref{shswgreenov}),  (\ref{deltaov}) and $G=E+\gtild$ we
obtain:
\begin{equation}
\begin{aligned} v(x,t)
&=\int^{+\infty}_{-\infty} \gtild_y(x,t;y)V_0(y)dy\\
&-\int^t_0
\int^{+\infty}_{-\infty}\gtild_y(x,t-s;y)\fcal(\fe, v,\ddp\delta)(y,s) dy \, ds,\\
&+ \int^t_0\int^{+\infty}_{-\infty}\gtild(x,t-s;y)\Psi(y,s) dy \,
ds.  \label{gtildov}
\end{aligned}
\end{equation}

 We are now in possession of the necessary tools
to state the following theorem:

\begin{theo} \label{thmmainover}
  Let  $(\hcal)$ and $(\dcal)$ hold, and
$|u_0|_{L^1}$, $|xu_0|_{L^1}$, $|u_0|_{L^\infty} \leq E_0$, $E_0$
sufficiently small (these assumption on $u_0$ are being inherited
by $v_0$). Assume the above setting and $u=v+\fe+\ddp\delta$, then
for any $\epsilon, \, 0<\epsilon<\frac18$,
\begin{equation}
|v(\cdot,t)|_{L^p}\le
 C E_0
(1+t)^{-\frac{1}{2}(1-1/p)-\frac 14},
\end{equation}
\begin{equation}
|\delta (t)|\le C E_0(1+t)^{-\frac 12+\epsilon},
\end{equation}
\begin{equation}
|\dot \delta (t)|\le C E_0(1+t)^{-1+\epsilon},
\end{equation}
 for any $p,$
$1 \leq p \leq \infty$, and with $C$ independent of $p$ (but
depending on $\epsilon$).
\end{theo}
An immediate consequence to this theorem is the following
corollary, which is almost (up to an $\epsilon$) Liu's result.
\begin{cor}\label{cor1}
\begin{equation}
|\utild - \ubar^{\delta_0}-\fe|_{L^p} \le \begin{cases}
(1+t)^{-\frac{1}{2}(1-1/p)-\frac 14} \quad &\textup{for}\, 1\le p
\le \frac{2}{1+8\epsilon},\\
(1+t)^{-\frac 12+\epsilon} \quad &\textup{for}\,  p \ge
\frac{2}{1+8\epsilon}.
\end{cases}
\end{equation}
\end{cor}
However, our approach yields more information about the behavior
of the perturbation, as we can track the shock location: by
(\ref{ishl}) and the comment right after, the following corollary
follows.
\begin{cor}\label{cor2}
\begin{equation}
|\utild - \ubar^{\delta_0 +\delta(t)} -\fe|_{L^p} \le
(1+t)^{-\frac{1}{2}(1-1/p)-\frac 14}
\end{equation}
for all $p$.
\end{cor}
\begin{proof}[Proof of Theorem \ref{thmmainover}]
Fixing $\epsilon,$ define
\begin{equation}
\begin{aligned}
     \zeta(t) := \sup_{0\le s \le t,1 \leq p \leq \infty}|u(\cdot,
    s)|_{L^p}(1+s)^{\frac12(1-\frac1p)+\frac14}
    &+ \sup_{0\le s \le t} |\delta (s)|(1+s)^{\frac 12-\epsilon }\\
    & + \sup_{0\le s
    \le t} |\dot\delta (s)|(1+s)^{1-\epsilon}.
\end{aligned} \label{zetadefov}
\end{equation}
Our aim is to show that
\begin{equation}
       \zeta (t) \leq C(E_0 + \zeta^2(t)) \label{zetaineq}
\end{equation}
and then use a straightforward continuous induction. Equivalent to
(\ref{zetaineq}) is
\begin{equation}
      |v(\cdot, s)|_{L^p} \leq C(E_0
       +\zeta^2(t)) (1+s)^{-\frac12(1-\frac1p)-\frac14}
\end{equation} and
\begin{equation}
    |\delta(s)| \leq C(E_0 +\zeta^2(t)) (1+s)^{-\frac 12+\epsilon }
\label{delbd}
\end{equation}
and a similar statement for $\dot\delta$. We need to take three
steps:

\textbf{step 1: bounds for $\mathbf{|v|_{L^p}}$:} By corollary
\ref{ubardecay}, we have $|\frac{\partial\ubar^\delta}{\partial
\delta_i}| \sim e^{-k|x|}$
and $|(\frac{\partial\ubar^\delta}{\partial \delta_i})_x| \sim
e^{-k|x|}$ for some $k>0$, hence, for $0\leq s\leq t$ we have:
\begin{eqnarray}
    |\ddp\delta(s)|_{L^p} \leq C \zeta(t) (1+s)^{-\frac12+\epsilon},\\
    |(\ddp)_x\delta(s)|_{L^p} \leq C \zeta(t) (1+s)^{-\frac12+\epsilon},\\
    |\fe(\cdot, s)|_{L^p} \leq C E_0 (1+s)^{-\frac12(1-\frac1p)},\\
    |\fe_x(\cdot, s)|_{L^p} \leq C E_0
    (1+s)^{-\frac12(1-\frac1p)-\frac12}.
\end{eqnarray}

Lemma \ref{wxlem} provides us with the necessary bounds for $v_x$.
Note that by lemma \ref{hKlem} and the above bounds, we have the
following bounds for $\mcal$ in (\ref{mcal}):
\begin{equation}
|\mcal (\fe, \ddp\delta)(\cdot,s)|_{L^p} \leq C (E_0 +
\zeta(t))(1+s)^{-\frac12(1-\frac1p)-\frac14}
\end{equation}
when $1\le p \le 2$. Therefore,
\begin{equation}
 |v_x(\cdot,s)|_{L^p} \le
 \begin{cases}
 C(E_0 +
\zeta(t))s^{-\frac12(1-\frac1p)-\frac14}, &\hbox{for}\ s \ge 1\\
 C(E_0 +
\zeta(t))s^{-\frac12}, & \hbox{for}\ s\le 1.\\
 \end{cases}\label{vxbounds}
 \end{equation}

Now let $\pstar =\frac{2}{1+ 8\epsilon},$ hence
$\frac12(1-\frac1\pstar)+\frac34=1-2\epsilon$. With  bounds for
$\fe, \delta, v$ and $v_x$ we obtain:
\begin{equation}
|\fcal (v, \fe, \ddp\delta)(\cdot, s)|_{L^p}\leq
C(E_0+\zeta(t)^2)s^{-\frac{1}{2}(1-1/p) - \frac 34}
\label{fcalbound}
\end{equation}
whenever $1\leq p \leq \pstar.$

 When $t \leq 1$, then
\begin{equation}
\begin{aligned} &|v(\cdot, t)|_{L^p}
\le C |v_0|_{L^p}\\
&+\int^t_0 |\gtild_y|_{L^1}|\fcal(v, \fe, \ddp\delta)|_{L^p}(s)ds\\
&+\int^t_0 (|\gtild|_{L^1}||\psi(y,s)|_{L^p}(s)ds\\
 &\le C E_0+(E_0+\zeta(t)^2) \int^t_0 (t-s)^{-\frac{1}{2}}
s^{-\frac{1}{2}}ds +C E_0 \int^t_0 (1+s)^{\frac 12}\\
&\le C (E_0 +\zeta(t)^2)\leq C (E_0 +\zeta(t)^2)(1+t)^{-\frac12(1-\frac1p)-\frac14}.\\
\end{aligned}
\end{equation}

 For $t \geq 1$ we use again Haussdorf-Young inequality
 to obtain:
\begin{equation}
\begin{aligned}
&|\int_{-\infty}^{+\infty} \gtild_y(x,t;y) V_0(y) dy|_{L^P}\\
&\leq  |V_0|_{L^1} |\gtild_y|_{L^p}\leq C E_0
t^{-\frac{1}{2}(1-1/p)-\frac 12},
\end{aligned} \label{b1}
\end{equation}
and
\begin{equation}
\begin{aligned}
|\int_{0}^{t/2} \int _{-\infty}^{+\infty}& \gtild_y(x,t-s;y)
\fcal(v,
\fe, \ddp\delta)(y,s) dy \, ds|_{L^p} \\
&\leq \int^{t/2}_0 |\gtild_y|_{L^p}|\fcal(v, \fe, \ddp\delta)(\cdot, s)|_{L^1}ds\\
&\leq C (E_0 + \zeta (t)^2)\int^{t/2}_0
(t-s)^{-\frac{1}{2}(1-1/p)-\frac{1}{2}}
s^{- \frac 34}ds\\
&\leq C(\zeta_0 +\zeta (t)^2)t^{-\frac{1}{2}(1-1/p) -\frac 14}\\
&\leq 2C(\zeta_0 +\zeta (t)^2)(1+t)^{-\frac{1}{2}(1-1/p) -\frac
14}.
\end{aligned}
\end{equation}
If $1\leq p \leq \pstar$, then
\begin{equation}
\begin{aligned}
|\int_{t/2}^t \int _{-\infty}^{+\infty}& \gtild_y(x,t-s;y)
\fcal(v,
\fe, \ddp\delta)(y,s) dy \, ds|_{L^p} \\
& \leq \int^t_{t/2}|\gtild_y|_{L^1}|\fcal(v, \fe, \ddp\delta)(\cdot, s)|_{L^p}ds\\
&\leq C (E_0 + \zeta (t)^2)\int^t_{t/2}(t-s)^{-\frac{1}{2}}
s^{-\frac{1}{2}(1-1/p) -
\frac 34} ds\\
&\leq C(E_0 +\zeta (t)^2)t^{-\frac{1}{2}(1-1/p) -\frac 14}\\
&\leq 2C(E_0 +\zeta (t)^2)(1+t)^{-\frac{1}{2}(1-1/p) -\frac 14}.
\end{aligned}
\end{equation}
If $p\geq \pstar,$ then choose $q$ so that $\frac1p+1=\frac1\pstar
+\frac1q$, and then
\begin{equation}
\begin{aligned}
|\int_{t/2}^t \int _{-\infty}^{+\infty}& \gtild_y(x,t-s;y)
\fcal(v,
\fe, \ddp\delta)(y,s) dy \, ds|_{L^p} \\
& \leq \int^t_{t/2}|\gtild_y|_{L^q}|\fcal(v, \fe, \ddp\delta)(\cdot, s)|_{L^{\pstar}}ds\\
&\leq C (E_0 + \zeta
(t)^2)\int^t_{t/2}(t-s)^{-\frac12(1-\frac1q)-\frac{1}{2}}
s^{-\frac{1}{2}(1-\frac1\pstar) -
\frac 34} ds\\
&\leq C(\zeta_0 +\zeta (t)^2)t^{-\frac{1}{2}(1-\frac1\pstar) -\frac 34+\frac{1}{2q}}\\
&= C(E_0 +\zeta (t)^2)t^{-\frac{1}{2}(1-1/p) -\frac 14}\\
&\leq 2C(E_0 +\zeta (t)^2)(1+t)^{-\frac{1}{2}(1-1/p) -\frac 14}.
\end{aligned}
\end{equation}
 It remains to show
$$
    |\int_0^t
\int_{-\infty}^{+\infty}\gtild(x,t-s;y) \Psi(y,s)dyds|_{L^p}\leq
E_0 (1+t)^{-\frac12(1-\frac1p)-\frac14}.$$

 By (\ref{ghatiov}),(\ref{ghatitarov}) and (\ref{ajabaaov})
 we have to estimate the following:
 First, the terms in the form $(\fe_i \fe_j \Gamma(x)(r_i^\pm,
r_j^\pm))_x,$ in (\ref{ghatiov}).  By lemma \ref{interaction},
$(\fe_i \fe_j \Gamma(x)(r_i^\pm, r_j^\pm))$ is of order $E^2_0
\bold{O} (e^{-\eta t}),$, so we can use similar calculations as
before to conclude that $|\int_{0}^t \int _{-\infty}^{+\infty}
\gtild_y(x,t-s;y) (\fe_i \fe_j \Gamma(x)(r_i^\pm, r_j^\pm)) dy \,
ds|_{L^p}=\bold O (E_0 (1+t)^{-\frac{1}{2}(1-1/p) -\frac 14})$

Next,  the terms in the second line of (\ref{ghatitarov}), i.e.,
in the form $((A(x)-A^-)\fe^i r_i^-)_x$ or similar forms, of which
lemma \ref{hKlem} together with (\ref{ABboundsov}) and similar
bounds for $\Gamma$ take care.

Finally, the terms in the form $\fe^i_{xx} b^-_{ij}r^-_j$ and
$(\fe^i)^2_{x} \Gamma^-_{jii}r^-_j,$  $i\neq j,$ in
(\ref{ajabaaov}). The $L_p$ norm of these terms is of order
$(1+t)^{-\frac12(1-\frac 1p)-\frac 12}$, so not decaying fast
enough to use calculations similar to what we have already done,
so we have to use the results in section \ref{S:preliminary}.
Assume $a_i^- <0$; we examine the integration of $(\fe^i)^2_{x}
r^-_j$ against the different terms of $S(x,t-s;y).$ As $l^{-tr}_k
r^-_j=0$ if $k\neq j$ the terms of concern in $S$ integrated
against $(\fe^i)^2_{xx} r^-_j$ are as following: For $y<0$, the
terms in the first line of (\ref{Sov}), gives us:
\begin{equation}
\int_{0}^t \int (4\pi \beta_j^-(t-s))^{-1/2}
e^{-(x-y-a_j^-(t-s))^2 / 4\beta_j^-(t-s)}(\fe_i)^2_ydy\,ds,
\label{gfgo}
\end{equation}
which, by (\ref{sect1cor}) and (\ref{phicor}), is of order
 $E_0(1+t)^{-\frac12(1-\frac1p)-\frac14}$ (as $i\neq j$). The second line of (\ref{Sov}) does not comprise
 anything. The third line of $S$ comprises:
\begin{equation}
 \int_{0}^t \int \chi_{\{t-s\ge 1\}}
[c^{j,-}_{k,-}] (4\pi \bar\beta_{jk}^- (t-s))^{-1/2}
e^{-(x-z_{jk}^-)^2 / 4\bar\beta_{jk}^- (t-s)} \left({\frac{e^{
-x}}{e^x+e^{-x}}}\right)(\fe_i)^2_y dy\,ds \label{chcoeff}
\end{equation}
with $a_k^->0, a_j^-<0, $ and $z_{jk}^-$ and $\bar\beta_{jk}^-$
computed at $t-s$. Note in this case the convection and diffusion
coefficients are not constant. To make a brief presentation of the
relevant calculations, we first notice that the biggest part in
(\ref{chcoeff}) is in the cone $a_j^-(t-s)/2\le x-z_{jk}^-\le
-a_j^-(t-s)/2,$ or equivalently, $3a_j^-(t-s)/2\le x- {a_j^-
y}/{a_k^-}\le a_j^-(t-s)/2$ (outside this cone we have a
negligible term). On this interval the diffusion coefficient
$\bar\beta_{jk}^-$ can be bounded from above by a constant
$\beta^*$.  The derivatives of $\bar\beta_{jk}^- (t-s)$ also can
be bounded from above similarly. As a consequence the $y$ and $t$
derivatives of this part of $S$ satisfy the bounds used in
proposition \ref{sect1main} in the cone just mentioned. Also we
make a change of coordinates $z= \frac{a_j^- y}{a_k^-}$ to see
that , this part of (\ref{chcoeff}) can be estimated the same way
one would estimate
$$\int\int g(x-z-a_j^-(t-s), \beta^*(t-s))g(z-a_i^-a_j^-s/a_k^-,
s)^2_ydyds.$$ using the same process as in proposition
\ref{sect1main} and subsequent results. Notice that
$a_i^-a_j^-/a_k^-\ne a_j^-$, i.e., the different speed of the
Gaussian kernel in the Green function and the diffusion wave.

 For $y>0$,  the corresponding terms in second line of $S$ in (\ref{Sov}) gives:
 \begin{equation}
 \int_{0}^t \int (4\pi \beta_k^+(t-s))^{-1/2}
e^{-(x-y-a_k^+(t-s))^2 / 4\beta_k^+(t-s)} \left({\frac
{e^{x}}{e^x+e^{-x}}}\right)(\fe_i)^2_y dy\,ds. \notag
\end{equation}
The case would be different from (\ref{gfgo}) if $a_k^+ = a_i^-$. In this case we have a term like:
 \begin{equation}
 \begin{aligned}
 \left({\frac {e^{x}}{e^x+e^{-x}}}\right)\int_{0}^t \int g(x-y-a(t-s),t-s)g^2(y-as,s) dy\,ds,\\
\end{aligned}
\end{equation}
with $a<0,$ which, by some elementary calculations,  is less than
or equal to:
\begin{equation}
 \begin{aligned}
 {\frac {Ce^{x}}{e^x+e^{-x}}}g(x-at, Mt).\\
\end{aligned}
\end{equation}
Now use lemma \ref{hKlem}.\\
 The terms in the remainder
$R(x,t;y)$ can be dealt with in a similar way.

\textbf{step 2: bounds for $\mathbf{\delta(t)}$:} In order to show
(\ref{delbd}) holds we investigate the integration of $e_{i_y}$
and $e_i$ against each term in $V_0 ,\fcal(v, \fe, \ddp\delta)$
and $\Psi,$ respectively. To that end, we will have the following
(assume $t\geq 1$).
\begin{equation}
\begin{aligned}
\int_{-\infty}^{+\infty}&e_{i_{y}}(y, t)V_0(y) dy\\
&\le |e_{i_{y}}|_{L^{\infty}}\,|V_0|_{L^1} dy\\
&\le Ct^{-\frac12},
\end{aligned}
\end{equation}
and
\begin{equation}
\begin{aligned}
|\int_0^t \int_{-\infty}^{+\infty}&e_{i_{y}}(y, t-s)\bold O (|v|^2)(y,s)dy\,ds| \\
&\leq C\int_0^{t/2} |e_{i_{y}}|_{L^{\infty}}(y,t-s)||v|^2|_{L^1}(y,s)ds\\
&\quad + C\int_{t/2}^t |e_{i_{y}}|_{L^1}(y,t-s)||v|^2|_{L^{\infty}}(y,s)ds\\
&\leq C\zeta^2(t)(\int_0^{t/2}(t-s)^{-\frac12}(1+s)^{-1}(y,s)ds
+ \int_{t/2}^t (1+s)^{-\frac 32}ds)\\
&\leq C\zeta^2(t)(1+t)^{-\frac12+\epsilon}.
\end{aligned}
\end{equation}
Similarly for $\bold O (|v||v_x|)$. For $\bold O((\ddp\delta)^2)$
we use lemma \ref{Kelem} to get
\begin{equation}
\begin{aligned}
|\int_0^t \int_{-\infty}^{+\infty}&e_{i_{y}}(y, t-s)\bold O (|\ddp\delta|^2)(y,s)dy\,ds \\
&\leq C \zeta^2(t)\int_0^t (1+s)^{-\frac12+\epsilon} \int_{-\infty}^{+\infty}e_{i_{y}}(y, t-s) (|\ddp|^2)(y,s)dy\,ds \\
&\leq C \zeta^2(t)\int_0^t (1+s)^{-\frac12+\epsilon} (t-s)^{-\frac12}e^{-\eta (t-s)}ds \\
&\leq C\zeta^2(t)(1+t)^{-\frac12}.
\end{aligned}
\end{equation}
For $\bold O (|\fe|^2)$ we use the fact that both $\fe$ and
$e_{i_{y}}$ are the summation of signals like convecting heat
kernels, moving away from shock. Hence using lemma \ref{KKKlem}
gives us:
\begin{equation}
\begin{aligned}
|\int_0^t \int_{-\infty}^{+\infty}e_{i_{y}}(y, t-s)\bold O
(|\fe|^2)(y,s)dy\,ds \leq C E_0 (1+t)^{-\frac12}.
\end{aligned}
\end{equation}
All the other terms in $\fcal(v,\fe,\ddp\delta)$ and $\Psi$ can be
treated with similar methods.

\textbf{step 3: bounds for $\mathbf{\dot\delta(t)}$:} Very similar
to the previous calculations for $\delta(t)$.
\end{proof}

\begin{rem}\label{liuinitial}
\textup{As for the conditions on initial data $u_0$, it is in fact
enough to assume (as Liu and others have done) that $u_0 = \bold O
(1+|x|)^{-\frac32}.$ To see briefly why this works notice that the
only part we should change in our argument is the linear part
(\ref{b1}). Now for this part, it is enough to consider the
convolution of a Gaussian signal, say g(x,t), against $v_0$. As
$|v_0|\sim (1+|x|)^{-\frac32}$ and $|V_0|\sim (1+|x|)^{-\frac12}$,
we consider $\int_{-\infty}^{+\infty}g_y(x-y, t)V_0(y)dy$  when
$|x|\le \sqrt{t},$ and $\int_{-\infty}^{+\infty}g(x-y, t)v_0(y)dy$
when $|x|\le \sqrt{t}.$ It is not difficult to observe, using
Howard's lemma \ref{holem}, that
$$\big|\int_{-\infty}^{+\infty}\gtild(x-y, t)v_0(y)dy\big|\le C\big(\chi_{\{|x|\le \sqrt{t}\}}
t^{-\frac12}(1+|x|)^{-\frac12} +  \chi_{\{|x|\ge
\sqrt{t}\}}(1+|x|)^{-\frac32}\big).$$ The necessary $L^1$ bounds
then follows immediately. To obtain $L^\infty$ bounds in the case
$|x|\le \sqrt{t}$ we use lemma \ref{stlem}, to see  that
\begin{equation}\notag\begin{aligned}
\int_{-\infty}^{+\infty}&g_y(x-y, t)V_0(y)dy\\
             &\le C t^{-\frac12}\int_{-\infty}^{+\infty}\tnim
             e^{-(x-y)^2/4t}y^{-\frac12}dy\\
             &\le C t^{-\frac34}\int_{-\infty}^{+\infty}
             e^{-y^2/8t}(\frac{y}{\sqrt{t}})^{-\frac12}\frac{dy}{\sqrt{t}}\\
             &\le Ct^{-\frac34}.
\end{aligned}\end{equation}
Other $L^p$ bounds follow using interpolation.}
\end{rem}
\begin{rem}\label{L2bas}\textup{ It is an easy observation that, in order to
have (\ref{fcalbound}) for $1\le p\le \pstar <2$, it is enough to
have $L^2$ bounds for $v_x$ in (\ref{vxbounds}), as $v_x$ is
always multiplied by a favorable term in $\fcal$. This will become
important when we consider the real viscosity case, as we have
only good $L^2$ bounds (using energy estimates) for $v_x$ in that
case . }
\end{rem}
\begin{rem}\textup{
The Analysis can go through in the case $df(u_\pm)$ is not
strictly hyperbolic, provided that $A_\pm$ and $B_\pm$ are
simultaneously symmetrizable, by replacing Green function bounds
with more general bounds given in proposition 5.10 of
\textbf{\cite{Z.5}}, and replacing the diffusion waves $\fe_i$'s
with the ``multi-mode diffusion waves" of \textbf{\cite{Ch}} and
\textbf{\cite{LZe}}.  The same remark is applicable in the real
viscosity case. }
\end{rem}

\section{Real viscosity case} \label{S:realvisc}

In this section we follow closely the notations and assumptions
used in \textbf{\cite{Z.5}}. Consider a general system of {\it
viscous conservation laws}
\begin{equation}
\begin{aligned}
 U_t+ F(U)_{x}&=
(B(U)U_{x})_{x}, \cr x\in \Bbb{R}; \,  &U, \, F\in \Bbb{R}^n; \,
B\in  \Bbb{R}^{n\times n}, \label{viscous}
\end{aligned}
\end{equation}

modeling flow in a compressible medium. We assume
\begin{equation}
U=\left(\begin{matrix}  u^I\\ u^{II}\end{matrix}\right), \quad
B=\left(\begin{matrix} 0 & 0\\
b_1 & b_2 \end{matrix}\right), \label{UB}
\end{equation}
$ u^I\in \BbbR^{n-r}$, $u^{II}\in \BbbR^r$,
and
\begin{equation}
Re \sigma  b_2 \ge \theta , \label{goodb}
\end{equation}
with $\theta>0$.

Again we consider the {\it  viscous shock wave} solutions of
(\ref{viscous}), which are in the form:
\begin{equation}
U(x,t)=\bar U(x), \quad \lim_{x\to \pm \infty} \bar U(x)=U_\pm,
\label{profile}
\end{equation}
satisfying the traveling-wave ordinary differential equation (ODE)
\begin{equation}
B(\bar U) \bar U'=F(\bar U)-F(U_-). \label{ODE}
\end{equation}
Considering the block structure of $B$, this can be written as:
\begin{equation}
F^I(u^I, u^{II})\equiv F^I(u_-^I, u_-^{II})\label{eq1}
\end{equation}
and
\begin{equation}
b_1(u^I)' + b_2(u^{II})'= F^{II}(u^I, u^{II}) - F^{II}(u_-^I,
u_-^{II}). \label{eq2}
\end{equation}


We assume that, by some invertible change of coordinates $U\to
W(U)$, possibly but not necessarily connected with a global convex
entropy, followed if necessary by multiplication on the left by a
nonsingular matrix function $S(W)$, equations (\ref{viscous}) may
be written in the {\it quasilinear, partially symmetric
hyperbolic-parabolic form}
\begin{equation}
\tilde A^0 W_t +  \tilde A W_{x}= (\tilde B W_{x})_{x} + G, \quad
W=\left(\begin{matrix} w^I
\\w^{II}\end{matrix}\right), \label{symm}
\end{equation}
 $w^I\in \BbbR^{n-r}$, $w^{II}\in \BbbR^r$, $x\in
\BbbR^d$, $t\in \BbbR$, where, defining $W_\pm:= W(U_\pm)$:\\
\medskip
(A1)\quad $\tilde A(W_\pm)$, $\tilde A_*:=\tilde A_{11}$, $\tilde
A^0$ are symmetric, $\tilde A^0 >0$.\\
\medskip
(A2)\quad No eigenvector of $ dF(U_\pm)$ lies in the kernel of $
B(U_\pm)$. (Equivalently, no eigenvector of $ \tilde A (\tilde
A^0)^{-1}(W_\pm)$ lies in the kernel of $\tilde B(W_\pm)$.)\\
\medskip
(A3) \quad $ \tilde B= \left(\begin{matrix} 0 & 0 \\ 0 & \tilde b
\end{matrix}\right) $, $ \tilde G= \left(\begin{matrix}  0 \\ \tilde g\end{matrix}\right) $,
with $ Re  \tilde b(W) \ge \theta $ for some $\theta>0$, for all
$W$, and $\tilde g(W_x,W_x)=\CalO(|W_x|^2)$.
\medskip
Here, the coefficients of (\ref{symm}) may be expressed in terms
of the original equation (\ref{viscous}), the coordinate change
$U\to W(U)$, and the approximate symmetrizer $S(W)$, as
\begin{equation}
\begin{aligned} \tilde A^0&:= S(W)(\partial U/\partial W),\quad
\tilde A:= S(W)d(\partial U/\partial W),\\
\tilde B&:= S(W)B(\partial U/\partial W), \quad G= -(dS W_{x})
B(\partial U/\partial W) W_{x}.
\end{aligned}
\label{coeffs}
\end{equation}
For examples about Navier--Stokes and Magnetohydrodynamic
equations, see \textbf{\cite{Z.4}}. Along with the above
structural assumptions, we make the technical
hypotheses:\\
\medskip
(H0) $F$, $B$, $W$, $S\in C^{s}$,
with $s\ge 5.$\\
\smallskip
(H1) The eigenvalues of $\ta_*$ are (i) distinct from $0$; (ii) of
common sign; and (iii) of constant multiplicity with respect to
$U$.\\
\smallskip
(H2) $\sigma (dF(U_\pm))$ real, distinct, and
nonzero.\\
\smallskip
(H3)  Local to $\bu(\cdot)$, solutions of
(\ref{profile})--(\ref{ODE}) form a smooth manifold $\{\bar
U^\delta(\cdot)\}$, $\delta \in \mathcal{U} \subset \BbbR^\ell$.
\medskip

Analogous to lemma \ref{mplem} and corollary \ref{ubardecay} we
have the following lemma proved in \textbf{\cite{MaZ.3}}.
\begin{lem}
Given (H1)--(H3), the endstates $U_\pm$ are hyperbolic rest points
of the ODE determined by (\ref{eq2}) on the $r$-dimensional
manifold (\ref{eq1}), i.e., the coefficients of the linearized
equations about $U^\pm$, written in local coordinates, have no
center subspace. In particular, under regularity (H0),
\begin{equation}
D_x^j D_\delta^i(\bar
U^\delta(x)-U_{\pm})=\bold{O}(e^{-\alpha|x|}), \quad \alpha>0, \,
0\le j\le 6, \,i=0,1,
\end{equation}
as $x\rightarrow\pm\infty$.
\end{lem}

We now recall some important ideas of Kawashima et al concerning
the smoothing effects of hyperbolic--parabolic coupling.  The
following results assert that hyperbolic effects can compensate
for degenerate viscosity $B$, as depicted by the existence of a
{\it compensating matrix} $K$.

\begin{lem} \label{skew} (\textbf{\cite{KSh}}) Assuming $A^0$, $A$, $B$
symmetric, $A^0>0$, and $B\ge 0$, the genuine coupling condition
\medskip
(GC)\quad  No eigenvector of $A$ lies in $\ker B$
\medskip
\noindent is equivalent to either of:\\
\medskip
(K1) \quad There exists a smooth skew-symmetric matrix function
$K(A^0,A,B)$ such that
\begin{equation}
\text{ Re }\left( K(A^0)^{-1}A + B \right)(U) >0. \label{skew}
\end{equation}\\
\medskip
(K2)  \quad For some $\theta>0$, there holds
\begin{equation}
\text{\rm Re }\sigma(-i\xi (A^0)^{-1}A-|\xi|^2 (A^0)^{-1}B) \le
-\theta |\xi|^2/(1+|\xi|^2), \label{symbol}
\end{equation} for all $\xi\in \mathbb{R}$.
\end{lem}

\begin{proof}
 These and other useful equivalent formulations are
established in \textbf{\cite{KSh}}; see also \textbf{\cite{ Z.5,
MaZ.4, Z.4}}.
\end{proof}

Now returning to the original equation (\ref{viscous}) with $\bu$
a shock solution we linearize around $ \bu$ exactly as we did in
section \ref{S:overc}, and we define again $A$ and $B$ in the same
manner:
\begin{equation}
B(x):= B(\bu(x)), \quad  A(x)V:= dF(\bu(x))V-dB(\bu(x))V\bu_x.
\label{AandBrv}
\end{equation}
Assume for  $A$ and $B$ the block structures:
$$A=\left(\begin{matrix}A_{11}\quad A_{12}\\A_{21}\quad A_{22}\end{matrix}\right),
B=\left(\begin{matrix}0& 0\\B_{21}& B_{22}\end{matrix}\right).$$
 The characteristics speeds $a_i^\pm,$
the left and right eigenvalues $l_i^\pm, r_i^\pm$ for $dF(u_\pm)$
and $\beta_i^\pm=l_i^\pm B^\pm r_i^\pm$ are all defined the same
way as before and again we have $\beta_i^\pm > 0.$

Also, let $a^{*}_j(x)$, $j=1,\dots,(n-r)$ denote the eigenvalues
of
$$
A_{*}:= A_{11}- A_{12} B_{22}^{-1}B_{21}, \label{A*}
$$
with $l^*_j(x)$, $r^*_j(x)\in \BbbR^{n-r}$ associated left and
right eigenvectors, normalized so that $l^{*t}_jr_j\equiv
\delta^j_k$. More generally, for an $m_j^*$-fold eigenvalue, we
choose $(n-r)\times m_j^* $ blocks $L_j^*$ and $R_j^*$ of
eigenvectors satisfying the {\it dynamical normalization}
$$
L_j^{*t}\partial_x R_j^{*}\equiv 0,
$$
along with the usual static normalization $L^{*t}_jR_j\equiv
\delta^j_kI_{m_j^*}$; as shown in Lemma 4.9,\textbf{
\cite{MaZ.1}}, this may always be achieved with bounded $L_j^*$,
$R_j^*$. Associated with $L_j^*$, $R_j^*$, define extended,
$n\times m_j^*$ blocks
$$
\mathcal{L}_j^*:=\left(\begin{matrix} L_j^* \\
0\end{matrix}\right), \quad \mathcal{R}_j^*:=
\left(\begin{matrix} R_j^*\\
-B_{22}^{-1}B_{21} R_j^*\end{matrix}\right). \label{CalLR}
$$
%
Eigenvalues $a_j^*$ and eigenmodes $\mathcal{L}_j^*$,
$\mathcal{R}_j^*$ correspond, respectively, to short-time
hyperbolic characteristic speeds and modes of propagation for the
reduced, hyperbolic part of degenerate system (\ref{viscous}).

Define local, $m_j\times m_j$ {\it dissipation coefficients}
$$
\eta_j^*(x):= -L_j^{*t} D_* R_j^* (x), \quad j=1,\dots,J\le n-r,
\label{eta}
$$
where
$$
\aligned
&{D_*}(x):= \\
& \, A_{12}B_{22}^{-1} \Big[A_{21}-A_{22} B_{22}^{-1} B_{21}+
A_{*} B_{22}^{-1} B_{21} + B_{22}\partial_x (B_{22}^{-1}
B_{21})\Big]
\endaligned
\label{D*}
$$
is an effective dissipation analogous to the effective diffusion
predicted by formal, Chapman--Enskog expansion in the (dual)
relaxation case.

At $x=\pm \infty$, these reduce to the corresponding quantities
identified by Zeng \textbf{\cite{Ze.1,LZe}} in her study by
Fourier transform techniques of decay to {\it constant solutions}
$(\bar u, \bar v) \equiv (u_\pm,v_\pm)$ of hyperbolic--parabolic
systems, i.e., of limiting equations
$$
U_t=L_\pm U:= -A_{\pm} U_x+ B_\pm U_{xx}. \label{limiting}
$$
As a consequence of dissipativity, (A2), we obtain (see, e.g.,
\textbf{\cite{Kaw, LZe, MaZ.3}}, or Lemma \ref{skew})
\begin{equation}
\beta_j^{\pm}>0, \quad Re \sigma(\eta_j^{*\pm})>0 \quad \text{\rm
for all $j$}. \label{goodbeta}
\end{equation}
However, note that the dynamical dissipation coefficient $D_*(x)$
{\it does not} agree with its static counterpart, possessing an
additional term $B_{22}\partial_x (B_{22}^{-1} B_{21})$, and so we
cannot conclude that (\ref{goodbeta}) holds everywhere along the
profile, but only at the endpoints. This is an important
difference in the variable-coefficient case; see Remarks 1.11-1.12
of \textbf{\cite{MaZ.3}} for further discussion.

We also make the following assumptions, necessary for linear
stability: \noindent Assumptions (D):

$(D1)$ \quad  $L$ has no ($L^2$, without loss of generality)
eigenvalues in $\{Re \lambda \ge 0\} \setminus \{0\}$.
\medskip

$(D2)$ \quad $\{r_j^\pm; a_j^\pm \gtrless 0\} \cup
\{\int_{-\infty}^{+\infty}\frac{\partial \bu^\delta}{\partial
\delta_i} dx; i=1, \cdots, \ell \}$ is a basis for $\mathbb{R}^n$,
with $\int_{-\infty}^{+\infty}\frac{\partial \bu^\delta}{\partial
\delta_i} dx$  computed at $\delta=0.$

\begin{prop} \label{greenbounds}\textbf{\cite{MaZ.3}}  Under assumptions (A1)--(A3),
(H0)--(H3), and (D1)--(D2),
 the Green distribution $G(x,t;y)$
associated with the linearized evolution equations  may be
decomposed as
$$
G(x,t;y)= H + E+  S + R, \label{ourdecomp}
$$
where, for $y\le 0$:
\begin{equation}
\begin{aligned} H(x,t;y)&:= \sum_{j=1}^{J} a_j^{*-1}(x) a_j^{*}(y)
\mathcal{R}_j^*(x) \zeta_j^*(y,t) \delta_{x-\bar a_j^* t}(-y)
\mathcal{L}_j^{*t}(y)\\
&= \sum_{j=1}^{J} \mathcal{R}_j^*(x) \mathcal{O}(e^{-\eta_0 t})
\delta_{x-\bar a_j^* t}(-y) \mathcal{L}_j^{*t}(y),
\end{aligned}
\label{multH}
\end{equation}
where the averaged convection rates $\bar a_j^*= \bar a_j^*(x,t)$
in (\ref{multH}) denote the time-averages over $[0,t]$ of
$a_j^*(x)$ along backward characteristic paths $z_j^*=z_j^*(x,t)$
defined by
$$
dz_j^*/dt= a_j^*(z_j^*), \quad z_j^*(t)=x, \label{char}
$$
and the dissipation matrix $\zeta_j^*=\zeta_j^*(x,t)\in
\BbbR^{m_j^*\times m_j^*}$ is defined by the {\it dissipative
flow}
$$
d\zeta_j^*/dt= -\eta_j^*(z_j^*)\zeta_j^*, \quad
\zeta_j^*(0)=I_{m_j}. \label{diss}
$$

 $E$ and $S$ have exactly the same form as in proposition \ref{greenovc},
and
\begin{equation}
\begin{aligned}
R(x,t;y)&=
\bold{O}(e^{-\eta(|x-y|+t)})\\
&+\sum_{k=1}^n
\bold{O} \left( (t+1)^{-1/2} e^{-\eta x^+}
+e^{-\eta|x|} \right)
t^{-1/2}e^{-(x-y-a_k^{-} t)^2/Mt} \\
&+
\sum_{a_k^{-} > 0, \, a_j^{-} < 0}
\chi_{\{ |a_k^{-} t|\ge |y| \}}
\bold{O} ((t+1)^{-1/2} t^{-1/2})
e^{-(x-a_j^{-}(t-|y/a_k^{-}|))^2/Mt}
e^{-\eta x^+}, \\
&+
\sum_{a_k^{-} > 0, \, a_j^{+}> 0}
\chi_{\{ |a_k^{-} t|\ge |y| \}}
\bold{O} ((t+1)^{-1/2} t^{-1/2})
e^{-(x-a_j^{+} (t-|y/a_k^{-}|))^2/Mt}
e^{-\eta x^-}, \\
\end{aligned}
\label{Rbounds}
\end{equation}
\begin{equation}
\begin{aligned}
R_y(x,t;y)&=
\sum_{j=1}^J \bold{O}(e^{-\eta t})\delta_{x-\bar a_j^* t}(-y)
+
\bold{O}(e^{-\eta(|x-y|+t)})\\
&+\sum_{k=1}^n
\bold{O} \left( (t+1)^{-1/2} e^{-\eta x^+}
+e^{-\eta|x|} \right)
t^{-1}
e^{-(x-y-a_k^{-} t)^2/Mt} \\
&+
\sum_{a_k^{-} > 0, \, a_j^{-} < 0}
\chi_{\{ |a_k^{-} t|\ge |y| \}}
\bold{O} ((t+1)^{-1/2} t^{-1})
e^{-(x-a_j^{-}(t-|y/a_k^{-}|))^2/Mt}
e^{-\eta x^+} \\
&+
\sum_{a_k^{-} > 0, \, a_j^{+} > 0}
\chi_{\{ |a_k^{-} t|\ge |y| \}}
\bold{O} ((t+1)^{-1/2} t^{-1})
e^{-(x-a_j^{+}(t-|y/a_k^{-}|))^2/Mt}
e^{-\eta x^-}, \\
\end{aligned}
\label{Rybounds}
\end{equation}
\begin{equation}
\begin{aligned}
R_x(x,t;y)&=
\sum_{j=1}^J \bold{O}(e^{-\eta t})\delta_{x-\bar a_j^* t}(-y)
+
\bold{O}(e^{-\eta(|x-y|+t)})\\
&+\sum_{k=1}^n
\bold {O} \left( (t+1)^{-1} e^{-\eta x^+}
+e^{-\eta|x|} \right)
t^{-1} (t+1)^{1/2}
e^{-(x-y-a_k^{-} t)^2/Mt} \\
&+
\sum_{a_k^{-} > 0, \, a_j^{-} < 0}
\chi_{\{ |a_k^{-} t|\ge |y| \}}
\bold{O}(t+1)^{-1/2} t^{-1})
e^{-(x-a_j^{-}(t-|y/a_k^-|))^2/Mt}
e^{-\eta x^+} \\
&+
\sum_{a_k^{-} > 0, \, a_j^{+} > 0}
\chi_{\{ |a_k^{-} t|\ge |y| \}}
\bold{O}(t+1)^{-1/2} t^{-1})
e^{-(x-a_j^{+}(t-|y/a_k^{-}|))^2/Mt}
e^{-\eta x^-}. \\
\end{aligned}
\label{Rxbounds}
\end{equation}
Moreover, for $|x-y|/t$ sufficiently large,
$|G|\le Ce^{-\eta t}e^{-|x-y|^2/Mt)}$ as in the strictly
parabolic case.
\end{prop}

Once again let $\gtild = S+R$ and define $e_i$'s as before.
Obviously the same bounds mentioned for $S$ and $e_i$'s in lemma
\ref{gebounds} hold here also. Furthermore we have the following
for $H$:

\begin{lem}
With the conditions in proposition \ref{greenbounds}, $H$
satisfies:
$$
|\int_{-\infty}^{+\infty} H(\cdot,t;y)f(y)dy|_{L^p} \le Ce^{-\eta
t} |f|_{L^p}, \label{Hbounds}
$$
$$
|\int_{-\infty}^{+\infty} H_x(\cdot,t;y)f(y)dy |_{L^p} \le
Ce^{-\eta t} |f|_{W^{1,p}}, \label{Hxbounds}
$$
for some $C, \eta > 0$, for any $p\ge 1$ and $f\in W^{1,p}$.
\end{lem}
\begin{proof}
See \textbf{\cite{MaZ.3, MaZ.4}}.
\end{proof}
From here on, almost everything would be very similar to the
strictly parabolic case in section \ref{S:overc}: by replacing
$\bu$ with $\bu^{\delta_0}$, for a small $\delta_0$, we may assume
that the initial perturbation has no mass at the $\D
\bu^\delta\over \D \delta_i$ directions, then we define diffusion
waves, $\fe_i$, exactly as in (\ref{ovfiithshock-}) and
(\ref{ovfiithshock+}), then $\fe =
\sum_{a_i^-<0}\fe_i+\sum_{a_i^+>0}\fe_i$, and once again
$V=\tu-\bu^{\delta_0}-\fe-\ddpu\delta(t)$ with $\delta$ to be
found (from now on we once again assume, without loss of
generality, $\delta_0=0$). The equalities (\ref{khatiov}) to
(\ref{ajabaaov}) are reproduced exactly as before, but with $u$
and $v$ replaced by $U$ and $V$, respectively. Furthermore, it is
easy to see that
\begin{equation}
\begin{aligned}
\fcal(V, \fe, \ddpu\delta(t))_x& =
\bold{O}\big(\fcal(V,\fe,\ddpu\delta)\\
&+
|(V+\fe+\ddpu\delta)_x||(v^{II}+\fe+\ddpu\delta)_x|\\
&+|V+\fe+\ddpu\delta||(v^{II}+\fe+\ddpu\delta)_{xx}|\big).
\label{fcalx}
\end{aligned}
\end{equation}
 The function $\delta(t)$ is defined as in
(\ref{deltaov}), and then, similar to (\ref{gtildov}), we have:
\begin{equation}
\begin{aligned} V(x,t)
&=\int^{+\infty}_{-\infty}(H+\gtild)(x,t;y)V(y,0)dy\\
&-\int^t_0
\int^{+\infty}_{-\infty}(H+\gtild)(x,t-s;y)\fcal(\fe, V,\ddp\delta)_y(y,s) dy \, ds,\\
&+ \int^t_0\int^{+\infty}_{-\infty}(H+\gtild)(x,t-s;y)\Psi(y,s) dy
\, ds. \label{gtildHrv}
\end{aligned}
\end{equation}

\begin{theo} \label{rvmain}
  Let  (A1)--(A3) and  (H0)--(H3), (D1)--(D2) hold,
and \linebreak $|U_0|_{L^1\cap L^\infty\cap H^3}$,
$|xU_0|_{L^1}\le E_0$, $E_0$ sufficiently small. Assume the above
setting and $U=V+\fe+\ddpu\delta$; then for any $\epsilon, \,
0<\epsilon<\frac18$,
\begin{equation}
|V(\cdot,t)|_{L^p}\le
 C E_0
(1+t)^{-\frac{1}{2}(1-1/p)-\frac 14},
\end{equation}
\begin{equation}
|\delta (t)|\le C E_0(1+t)^{-\frac 12+\epsilon},
\end{equation}
\begin{equation}
|\dot \delta (t)|\le C E_0(1+t)^{-1+\epsilon},
\end{equation}
 for any $p,$
$1 \leq p \leq \infty$, and with $C$ independent of $p$ (but
depending on $\epsilon$).
\end{theo}
\begin{rem}\textup{Corollaries similar to \ref{cor1} and \ref{cor2}
are valid here also.}
\end{rem}
\begin{proof}
Fixing $\epsilon,$ define
\begin{equation}
\begin{aligned}
     \zeta(t) := \sup_{0\le s \le t,1 \leq p \leq \infty}|V(\cdot,
    s)|_{L^p}(1+s)^{\frac12(1-\frac1p)+\frac14}
    + &\sup_{0\le s \le t} |\delta (s)|(1+s)^{\frac 12-\epsilon } \\
    + &\sup_{0\le s
    \le t} |\dot\delta (s)|(1+s)^{1-\epsilon}.
\end{aligned} \label{zetadefrv}
\end{equation}
To show
\begin{equation}
       \zeta (t) \leq C(E_0 + \zeta^2(t)) \label{zetaineqrv}
\end{equation}
we need to show
\begin{equation}
      |V(\cdot, s)|_{L^p} \leq C(E_0
       +\zeta^2(t)) (1+s)^{-\frac12(1-\frac1p)-\frac14}
\end{equation} and
\begin{equation}
    |\delta(s)| \leq C(E_0 +\zeta^2(t)) (1+s)^{-\frac 12+\epsilon
    }
\label{delbd}
\end{equation}
and a similar statement for $\dot\delta$. From here on the proof
goes very much similarly to the proof of theorem
\ref{thmmainover}, except for two issues: first, here we do not
have a lemma similar to lemma \ref{wxlem}, as the short time
estimates there need the strict parabolic hypothesis. Instead we
have to use some energy estimates in order to control the
derivatives of $V$. Using this method we will find out that, under
the assumptions of the problem,
\begin{equation}|V(\cdot, t)|_{H^3} \le C(E_0
+\zeta(t))(1+s)^{-\frac12}. \label{vhnorm}
\end{equation}
This in turn will implies (\ref{fcalbound}), for $1\le p \le
\pstar,$ with $\pstar$ as before (see remark \ref{L2bas}).
 The other difference is that, here we have
  the extra term
$H$ in the Green function decomposition, but we do not have any
bounds for $H_y$. Hence we have to compute
$$\int_0^t\int_{-\infty}^{+\infty}H(x,t;y)\fcal_y(y,s)dyds.$$
By (\ref{fcalx}) and (\ref{vhnorm}), we obtain
$$|\fcal_x(\cdot, s)|_{L^p} \le C(E_0
+\zeta(t)^2) s^{-\frac12(1-\frac1p)-\frac34},$$ $1\le p \le
\pstar.$ This with lemma \ref{Hbounds} provides us with necessary
bounds.

It remains to show that (\ref{vhnorm}) holds.  Let $\tilde U-\bar
U = V+\fe+\ddpu\delta$, and $W := \tilde W-\bar
W-\feh-\ddpw\delta$, with $\feh = d\bw \fe$. Notice also that
$d\bw(\ddpu)=\ddpw$ (we use the notation $d\bw=dW(\bu)$,
$d\bu=dU(\bw)$, etc).

\medskip
 \textbf{Claim:} $|V|_{H^r}\sim |W|_{H^r}.$\\
\emph{proof of the claim:} Note that $\tu -\bu = \duav(\tw-\bw),$
where $\duav = \int_0^1 dU(\bw +\theta(\tw-\bw)) d\theta$. Now
using the facts that $\ddpu=d\bu \ddpw$ and $\fe=d\bu \feh$ we
deduce:
$$V = \duav(W) +(\duav-d\bu)\feh + (\duav-d\bu)\ddpw\delta.$$
This, with a similar argument in the reverse direction proves our
claim.

From here on we follow closely the argument presented by Zumbrun
in \textbf{\cite{Z.5}} (see also \textbf{\cite{MaZ.4}}), with some
necessary modification to handle our more complicated case, e.g.,
key cancelations  in (\ref{D}) and (\ref{quasipert2}) and the term
$\xi$ in (\ref{whbounds}), which has no counterpart in
\textbf{\cite{Z.5}}. First, we introduce the weighted norms and
inner product
 \begin{equation}
|f|_\alpha:= |\alpha^{1/2}f|_{L^2}, \quad |f|_{H^s_\alpha}:=
\sum_{r=0}^s |\partial_x^r f|_\alpha, \quad \langle
f,g\rangle_\alpha:= \langle \alpha f,g\rangle_{L^2},
\label{alphanorm}
\end{equation}
$\alpha(x)$ scalar, uniformly positive, and uniformly bounded. For
the remainder of this section, we shall for notational convenience
omit the subscript $\alpha$, referring always to $\alpha$-norms or
-inner products unless otherwise specified. For later reference,
we note the commutator relation
 \begin{equation}
\langle f,g_x\rangle= -\langle f_x+ (\alpha_x/\alpha)f,g\rangle,
\label{commutator}
 \end{equation}
and the related identities
 \begin{equation}
\langle f,Sf_x\rangle= - (1/2)\langle f,\big(S_x
+(\alpha_x/\alpha) S\big)f\rangle, \label{symm1}
 \end{equation}
 \begin{equation}
\langle f,(Sf)_x\rangle=
 (1/2)\langle f,\big(S_x
-(\alpha_x/\alpha) S\big)f\rangle, \label{symm2}
 \end{equation}

valid for symmetric operators $S$.

By (H1)(ii), we have that $\bar A_{11}(\bar A^0_{11})^{-1}$ has
real spectrum of uniform sign, without loss of generality
negative, so that the similar matrix
$$
(\bar A^0_{11})^{-1/2}\bar A_{11}(\bar A^0_{11})^{-1/2}= (\bar
A^0_{11})^{-1/2} \bar A_{11}( \bar A^0_{11})^{-1} (\bar
A^0_{11})^{1/2}
$$
has real, negative spectrum as well.  (Recall, $\bar A^0_{11}$ is
symmetric negative definite as a principal minor of the symmetric
negative definite matrix $\bar A^0$.) It follows that $\bar
A_{11}$ itself is uniformly symmetric negative definite, i.e.,
 \begin{equation}
\bar A_{11}\le  -\theta<0. \label{transverse}
 \end{equation}

Defining $\alpha$, following Goodman \textbf{\cite{Go}}, by the
ODE
 \begin{equation}
\alpha_x= C_*|\bar W_x|\alpha, \quad \alpha(0)=1, \label{alphaODE}
 \end{equation}

where $C_*>0$ is a large constant to be chosen later, we have by
(\ref{transverse})
 \begin{equation}
(\alpha_x/\alpha) \bar A_{11} \le -C_*\theta |\bar W_x|.
\label{goodterm}
 \end{equation}
Note, because $|\bar W_x|\le Ce^{-\theta|x|}$, that $\alpha$ is
indeed positive and bounded from both zero and infinity, as the
solution of the simple scalar exponential growth equation
(\ref{alphaODE}).

\textbf{Energy estimates for $W$:}
\begin{equation}
\begin{aligned}
\tilde A^0 W_t + \tilde A W_x - (\tilde B& W_x)_x\\
= &-(\ta-\ba)\bw_x + ((\tb - \bb)\bw_x)_x \\
&- \ba \ddpwx\delta +\Big(\bb\ddpwx\Big)_x\delta\\
&-(\ta -\ba)\ddpwx\delta +\Big((\tb-\bb)\ddpwx\Big)_x\delta\\
&-\tao\ddpw\dot\delta -\tao \feh_t-\ta\feh_x+(\tb\feh_x)_x
\end{aligned}
\label{quasipert}
\end{equation}
 where
\begin{equation}
\begin{split}
\tilde A^0:=A^0(\tilde W),  \quad \tilde A:=A(\tilde W), \quad
\tilde B:=B(\tilde W);\\  \bar A^0:=A^0(\bar W), \quad \bar
A:=A(\bar W), \quad \bar B:=B(\bar W);
\end{split} \label{tildebar}
\end{equation}
(notice that we dropped tilde signs from $\tilde A^0, \tilde A$
and $\tilde B$ in (\ref{coeffs}) and, with a slight abuse of
notation, used them here differently). We want to write the right
hand side of (\ref{quasipert}) in the form $\emi+(\emii)_x +\xi$,
where $\emi$ and $\emii$ depending on $W$ and "behaving well
enough", and $\xi$ is a remainder decaying fast enough.

Beginning from the last line in (\ref{quasipert}), we write
$$\tao\ddpw\dot\delta = (\tao-
A^0(\bw+\feh+\ddpw\delta))\ddpw\dot\delta +
A^0(\bw+\feh+\ddpw\delta)\ddpw\dot\delta$$

Now $A^0(\bw+\feh+\ddpw)\ddpw\dot\delta$ goes into $\xi$ and
$(\tao- A^0(\bw+\feh+\ddpw))\ddpw\dot\delta$ goes into $\emi$.
Notice that
$$\tao-
A^0(\bw+\feh+\ddpw\delta)= \int_0^1 d A^0
(\bw+\feh+\ddpw\delta+\theta W)d\theta W$$ In a similar fashion
$\tao \feh_t$ and $\ta\feh_x$ each gives rise to a term which goes
into $\emi$ and the other which goes into $\xi$. Similarly
$(\tb\feh_x)_x$ comprises two terms, one of which is absorbed by
$\emii$ and the other by $\xi.$

Using \begin{equation}d\ba \ddpw\bw_x +\ba\ddpwx=(d\bb\ddpw
\bw_x)_x +(\bb\ddpwx)_x\label{D}\end{equation}
 we can write the second and third lines of
(\ref{quasipert}) in the form:
\begin{equation}
\begin{aligned}
(\ta-\ba)\bw_x - &((\tb - \bb)\bw_x)_x + \ba \ddpwx\delta -(\bb\ddpwx)_x\delta\\
= (&\ta-A(\bw+\feh+\ddpw\delta)) \bw_x \\
+&\left(A(\bw+\feh+\ddpw\delta)-\ba-d\ba(\feh+\ddpw\delta)\right)\bw_x + d\ba\feh\bw_x\\
-&\left((\tb-B(\bw+\feh+\ddpw\delta)) \bw_x\right)_x\\
+&\left((B(\bw+\feh+\ddpw\delta)-\bb-d\bb(\feh+\ddpw\delta))\bw_x\right)_x \\
-&(d\bb\feh\bw_x)_x\\
\end{aligned}
\label{quasipert2}
\end{equation}
the term in the second line of (\ref{quasipert2}) goes into
$\emi$, the third line and the fifth lines go into $\xi$, and the
fourth line goes into $\emii.$ The fourth line of
(\ref{quasipert}) can be dealt with in a similar way.

To summarize, we were able to write equation (\ref{quasipert}) in
the form:
\begin{equation}
\tilde A^0 W_t + \tilde A W_x - (\tilde B
W_x)_x=\emi+(\emii)_x+\xi(x,t) \label{leibnitz}
\end{equation}
where $\emi, \emii$ are dependent on $W$;
\begin{equation}
\begin{aligned}
\emi =& -(\ta^0- A^0(\bw+\feh+\ddpw\delta))\ddpw\dot\delta\\
      &-(\ta^0-A^0(\bw+\feh+\ddpw\delta))\feh_t\\
      &-(\ta-A(\bw+\feh+\ddpw\delta))\feh_x\\
      &-(\ta-A(\bw+\feh+\ddpw\delta)) \bw_x \\
      &-(\ta-A(\bw+\feh+\ddpw\delta))\ddpwx\delta.
\end{aligned}\label{emidef}
\end{equation}

\begin{equation}
\begin{aligned}
\emii =&(\tb-B(\bw+\feh+\ddpw\delta))\feh_x\\
       &+ (\tb-B(\bw+\feh+\ddpw\delta)) \bw_x\\
       &+(\tb-B(\bw+\feh+\ddpw\delta))\ddpwx\delta.
\end{aligned}\label{emiidef}
\end{equation}
Using $\tilde A= \bar A + \CalO(\zeta)$, $\tilde A_x= \CalO(|\bar
W_x| + \zeta)$ we can see:
\begin{equation}
|\langle \D_x^r W, \D_x^r\emi \rangle| \le C\langle  \D_x^r W,
\bw_x \D_x^r W\rangle + C\zeta
|W|^2_{H^r}+C|W|^2_{H^{r-1}}\label{emi}
\end{equation}
 for $r=0,\cdots, 3$. $\emii$ has the block
form: $\left(\begin{matrix}
0&0\\0&\mathcal{M}_{22}\end{matrix}\right) $, hence using Young's
inequality and the block structure of $\emii$,
\begin{equation}
|\langle \D_x^r W, \D_x^r(\emii)_x \rangle| \le C
\mu^{-1}|w^{II}|_{H^r}^2+ \mu |w^{II}|_{H^{r+1}}^2 \label{emii}
\end{equation}
for $\mu$ arbitrarily small. $\xi(x,t)$ is independent of $W$,
\begin{equation}
\begin{aligned}
\xi=&- A^0(\bw+\feh+\ddpw\delta)\ddpw\dot\delta\\
      &-A^0(\bw+\feh+\ddpw\delta)\feh_t\\
      &-A(\bw+\feh+\ddpw\delta)\feh_x\\
      &+(B(\bw+\feh+\ddpw\delta)\feh_x)_x\\
      &-\left(A(\bw+\feh+\ddpw\delta)-\ba-d\ba(\feh+\ddpw\delta)\right)\bw_x + d\ba\feh\bw_x\\
      &-\left((B(\bw+\feh+\ddpw\delta)-\bb-d\bb(\feh+\ddpw\delta))\bw_x\right)_x \\
      &+(d\bb\feh\bw_x)_x\\
      &-(A(\bw+\feh+\ddpw\delta)-\ba)\ddpwx\delta\\
      &+((B(\bw+\feh+\ddpw\delta)-\ba)\ddpwx)_x\delta.
\end{aligned}\label{xidef}
\end{equation}
 $\xi$
has the good property that, with differentiation with respect to
$x$, its rate of decay remains the same. Hence we transfer all the
derivatives to $\xi$:
\begin{equation}
|\langle D_x^r W,  D_x^r\xi \rangle| \le C (|W|_{L^2}^2 +
|\xi|_{H^{2r}}^2)\label{xi}
\end{equation}
Notice that under the assumptions and definitions of Theorem
\ref{rvmain},
\begin{equation}
|\xi(\cdot, s)|_{H^r} \le C(E_0+\zeta(t))(1+s)^{-\frac34}
\end{equation}
for any $r$ and $0\le s \le t$. Now the following lemma provides
us with necessary bounds we need for $|W|_{H^3}$:
\begin{lem}
Under the hypotheses of Theorem \ref{rvmain} let $W_0\in H^3$, and
suppose that, for $0\le t\le T$, both the supremum of $
|\dot\delta|, |\delta|$ and the $W^{2,\infty}$ norm of the
solution $W=(w^I,w^{II})^t$  remain bounded by a sufficiently
small constant $\zeta>0$. Then, for all $0\le t\le T$,
\begin{equation}
|W(t)|_{H^3}^2\le C |W(0)|^2_{H^3}e^{-\theta t} +C\int_0^t
e^{-\theta_2 (t-\tau )}(|W|_{L^2}^2+ |\xi|_{H^6}^2)(\tau)\, d\tau.
\label{whbounds}
\end{equation}
\end{lem}
We first carry out a complete proof in the more straightforward
case that the equations may be globally symmetrized , i.e., with
conditions (A1)--(A3) replaced by the following global versions,
indicating afterward by a few remarks the changes needed to carry
out the proof in the general
case.\\
\medskip
(A1')\quad $\tilde A^j$, $\tilde A^{jk}_*:=\tilde A^{jk}_{11}$,
$\tilde A^0$ are symmetric, $\tilde A^0 >0$.
\\
\medskip
(A2')\quad No eigenvector of $\sum \xi_j dF^j (U)$ lies in the
kernel of $\sum \xi_j \xi_k B^{jk}(U)$, for all nonzero $\xi\in
\BbbR^d$. (Equivalently, no eigenvector of $\sum \xi_j \tilde A^j
(\tilde A^0)^{-1}(W)$ lies in the kernel of $\sum \xi_j \xi_k
\tilde B^{jk}(W)$.)\\
\medskip
(A3') \quad
\begin{equation}
\tilde B^{jk}= \left(\begin{matrix} 0 & 0 \\
0 & \tilde b^{jk}
\end{matrix}\right), \end{equation}
with $ Re \sum \xi_j \xi_k \tilde{b}^{jk}(W) \ge \theta |\xi|^2$
for some $\theta>0$, for all $W$ and all $\xi\in \BbbR^d$, and
$\tilde G\equiv 0$.\\
\medskip
 To prove (\ref{whbounds}), we  carry out a series of successively higher order
energy estimates of the type formalized by Kawashima
\textbf{\cite{Kaw}} and used extensively by K. Zumbrun et al (see
\textbf{\cite{Z.5}}, also \textbf{\cite{MaZ.3, MaZ.4}}) The origin
of this approach goes back to \textbf{\cite{Kan,MNi}} in the
context of gas dynamics; see, e.g., \textbf{\cite{HoZ.1}} for
further discussion/references.

Let $\tilde K$ denote the skew-symmetric matrix described in Lemma
\ref{skew} associated with $\tilde A^0$, $\tilde A$, $\tilde B$,
satisfying
$$
\tilde K(\tilde A^0)^{-1}\tilde A+\tilde B>0.
$$
Then, regarding $\tilde A^0$, $\tilde K$, we have
\begin{equation}
\begin{aligned} \tilde A^0_x&= dA^0(\tilde W) \tilde W_x, \,\, \tilde
K_x= dK(\tilde W) \tilde W_x, \,\, \tilde A_x= dA(\tilde W) \tilde
W_x, \quad \tilde B_x= dB(\tilde W) \tilde W_x,
\\
\tilde A^0_t&= dA^0(\tilde W) \tilde W_t , \,\, \tilde K_t=
dK(\tilde W) \tilde W_t, \quad \tilde A_t= dA(\tilde W) \tilde
W_t, \quad \tilde B_x= dB(\tilde W) \tilde W_t. \\
\label{AKbounds}
\end{aligned}
\end{equation}

Now:
 \begin{equation} |\tilde
W_x|=|W_x + \bar W_x+\feh_x+\ddpwx\delta|\le |W_x|+ |\bar
W_x|+|\feh_x|+|\ddpwx\delta|
 \label{xbound}
\end{equation}

and, from the equation for $\ta$,
\begin{equation}
\begin{aligned} |\tilde W_t|&
\le C( |\tilde W_x| + |\tilde w^{II}_{xx}|).
\end{aligned}
\label{tbound} \end{equation}

Thus,  in particular it follows that

\begin{equation}
\begin{aligned} |\dot\delta|, \, |\tilde A^0_x|,\, |\tilde A^0_{xx}|,\,
|\tilde K_x|,\, |\tilde K_{xx}|,\, |\tilde A_x|,\, |\tilde
A_{xx}|,\, |\tilde B_x| ,\, |\tilde B_{xx}| ,\, |\tilde A^0_t|,\,
&|\tilde K_t|,\, |\tilde A_t|,\,
|\tilde B_t|  \\
&\le C(\zeta+|\bar U_x|).
\end{aligned}
\label{smallness} \end{equation}

 In what follows, we
shall need to keep careful track of the distinguished constant
$C_*$.

Computing
 \begin{equation}
-\langle W, \tilde AW_x\rangle= (1/2)\langle W, (\tilde A_x +
(\alpha_x/\alpha)\tilde A) W\rangle \label{one}
 \end{equation}

and expanding $\tilde A= \bar A + \CalO(\zeta)$, $\tilde A_x=
\CalO(|\bar W_x| + \zeta)$, we obtain by (\ref{goodterm}) the key
property
\begin{equation}
\begin{aligned} -\langle W, \tilde AW_x\rangle&=
(1/2)\langle w^I, (\alpha_x/\alpha) \bar A_{11} w^I\rangle\\
&\qquad +\CalO\Big( \langle |\bar W_x||W|,|W|\rangle +\langle
(\alpha_x/\alpha) |W|,\zeta |W|+ |w^{II}|\rangle
\Big)\\
&\le -(C_*\theta/3)\langle |w^I|, |\bar W_x| |w^I|\rangle + C\zeta
|w^I|^2 + C(C_*)|w^{II}|^2,
\end{aligned}
\label{goodtransprelim}
\end{equation}

  by which we shall control transverse
modes, provided $C_*$ is chosen sufficiently large, or, more
generally,
 \begin{equation}
-\langle \partial_x^kW, \tilde A\partial_x^kW_x\rangle \le
-(C_*\theta/3)\langle |\partial_x^k w^I|, |\bar W_x| |\partial_x^k
w^I|\rangle + C\zeta |\partial_x^k w^I|^2 + C(C_*)|\partial_x^k
w^{II}|^2. \label{goodtrans}
 \end{equation}

Here and below, $C(C_*)$ denotes a suitably large constant
depending on $C_*$, while $C$ denotes a fixed constant independent
of $C_*$: likewise, $\CalO(\cdot)$ indicates a bound independent
of $C_*$.

\medskip

{\bf Zeroth order ``Friedrichs-type'' estimate.} We first perform
a standard, zeroth- and first-order ``Friedrichs-type'' estimate
for symmetrizable hyperbolic systems \textbf{\cite{Fri}}. Taking
the $\alpha$-inner product of $W$ against (\ref{leibnitz}), we
obtain after rearrangement, integration by parts using
(\ref{commutator})--(\ref{symm1}), and several applications of
Young's inequality, the energy estimate
\begin{equation}
\begin{aligned} \frac{1}{2}\langle  W_{},\tilde A^0 W_{}\rangle_t &=
\langle W_{},\tilde A^0 W_{t}\rangle +\frac{1}{2}\langle
W_{},\tilde A^0_t W_{}\rangle
\\
&= -\langle W,\tilde A W_x\rangle
  +
\langle W,(\tilde B W_x )_x \rangle +\langle W, \emi \rangle
+ \langle W, (\emii)_x\rangle\\
&\qquad+\langle W,\xi\rangle
+\frac{1}{2}\langle  W_{},\tilde A^0_t W_{}\rangle\\
&= \frac{1}{2}\langle W,(\tilde A_x +(\alpha_x/\alpha)\tilde A)
W\rangle
  - \langle W_x- (\alpha_x/\alpha)W,\tilde B W_x  \rangle
\\& \qquad+\langle W, \emi \rangle - \langle W, \emii\rangle
+\langle W,\xi\rangle
+\frac{1}{2}\langle  W_{},\tilde A^0_t W_{}\rangle\\
&\le - \langle W_{x},\tilde BW_{x}\rangle
+ C(C_*)\int \alpha \Big( (|W_x|+|\bar W_x|)|W|^2\\
&\qquad + |w^{II}_x||W|(|W_x|+ |\bar W_x|) + C|W|_{L^2}^2 +
\mu|w^{II}|_{H^1}^2+
 |\langle W, \xi\rangle| \\
&\le -  \theta|w^{II}|_{H^1}^2 + C(C_*)\left( |W|^2_{L^2}
+|\xi|_{L^2}^2 \right). \\
\end{aligned}
\label{sym} \end{equation}
Here, we used boundedness of $|\D_x^r \bw|$ and $|W_x|$ and also
the inequalities (\ref{emi}), (\ref{emii}) and (\ref{xi}).

\smallskip

{\bf First order ``Friedrichs-type'' estimate.} For first and
higher derivative estimates, it is crucial to make use of the
favorable terms (\ref{goodtrans}) afforded by the introduction of
$\alpha$-weighted norms. Differentiating (\ref{leibnitz}) with
respect to $x$, taking the $\alpha$-inner product of $W_x$ against
the resulting equation, and substituting the result into the first
term on the righthand side of
 \begin{equation}
\frac{1}{2}\langle W_{x},\tilde A^0 W_{x}\rangle_t = \langle
W_{x},(\tilde A^0 W_{t})_x\rangle -\langle W_{x},\tilde A^0_x
W_{t}\rangle +\frac{1}{2}\langle  W_{x},\tilde A^0_t W_{x}\rangle,
\label{one}
 \end{equation}
we obtain after various simplifications and integrations by parts:
\begin{equation}
\begin{aligned} \frac{1}{2}\langle W_{x},\tilde A^0 W_{x}\rangle_t &=
-\langle W_x,(\tilde A W_x)_x\rangle
  +
\langle W_x,(\tilde B W_x )_{xx} \rangle
+\langle W_x, (\emi)_x \rangle \\
&\quad+ \langle W_x, (\emii)_{xx}\rangle
+\langle W_x, \xi_x\rangle\\
&\quad -\langle W_{x},\tilde A^0_x W_{t}\rangle
+\frac{1}{2}\langle  W_{x},\tilde A^0_t W_{x}\rangle\\
&= -\langle W_x,\tilde A W_{xx}\rangle -\langle W_x,\tilde A_x
W_{x}\rangle\\
&\qquad - \langle W_{xx}+(\alpha_x/\alpha)W_x,\tilde B W_{xx}
+ \tilde B_x W_{x}\rangle\\
&\qquad+\langle W_x, (\emi)_x \rangle
 - \langle W_{x}+ (\emii)_{xx}\rangle+\langle W_x, \xi_x\rangle\\
&\qquad   -\langle W_{x},\tilde A^0_x W_{t}\rangle
+\frac{1}{2}\langle  W_{x},\tilde A^0_t W_{x}\rangle.\\
\end{aligned}
\label{nextsym0} \end{equation}
Estimating the first term on the righthand side of
(\ref{nextsym0}) using (\ref{goodtrans}), $k=1$, and substituting
$(\tilde A^0)^{-1}$ times (\ref{leibnitz}) into the second to last
term on the righthand side of (\ref{nextsym0}), we obtain by
(\ref{smallness}) plus various applications of Young's inequality
the next-order energy estimate:
\begin{equation}
\begin{aligned} \frac{1}{2}\langle W_{x},\tilde A^0 W_{x}\rangle_t &\le
-\langle W_x,\tilde A W_{xx}\rangle
- \langle W_{xx},\tilde BW_{xx}\rangle \\
&\qquad
+C(C_*)\langle |W^{II}_x|+ \zeta |W_x|,(|W|+|W_x|)|\bar W_x|+|w^{II}_{xx}|\rangle\\
&\qquad+ C \langle (|W_x|+|w^{II}_{xx}|),|\bar W_x|(|W|+|W_x|) \rangle\\
&\qquad +C\langle W_x, \bw_x W_x\rangle +C\zeta|W_x|^2_{L^2}\\
&\qquad +C\mu^{-1}
|w^{II}_x|^2_{L^2}+\mu|w^{II}_{xx}|^2_{L^2}+ C(C_*)(|W|_{L^2}+|\xi|^2_{H^2})\\
\\
&\le -(\theta/2)|w^{II}_{xx}|_{L^2}^2 - (C_*\theta/4) \langle
|\bar W_x||w^I_x|,|w^I_x|\rangle + C(C_*)\zeta |w^I_x|_{L^2}^2\\
&\qquad + C(C_*)|w^{II}_x|_{L^2}^2
+ C(C_*)(|W|_{L^2}+|\xi|^2_{H^2}) ,\\
\end{aligned}
\label{nextsym} \end{equation}
 provided $C_*$ is sufficiently large and $\zeta, \mu$
sufficiently small.

\smallskip

{\bf First order ``Kawashima-type'' estimate.} Next, we perform a
``Kawashima-type'' derivative estimate. Taking the $\alpha$-inner
product of $W_{x}$ against $\tilde K(\tilde A^0)^{-1} $ times
(\ref{leibnitz}), and noting that (integrating by parts, and using
skew-symmetry of $\tilde K$)
\begin{equation}
\begin{aligned} \frac{1}{2}\langle W_{x},\tilde KW_{} \rangle_t &=
 \frac{1}{2}\langle W_{x},\tilde KW_{t} \rangle+
\frac{1}{2}\langle W_{xt},\tilde KW_{} \rangle
+\frac{1}{2}\langle W_{x},\tilde K_tW_{} \rangle\\
&= \frac{1}{2}\langle W_{x},\tilde KW_{t} \rangle
-\frac{1}{2}\langle W_{t},\tilde KW_{x} \rangle\\
&\quad -\frac{1}{2}\langle  W_{t},\big(\tilde
K_x+(\alpha_x/\alpha)\big) W \rangle
+\frac{1}{2}\langle W_{x},\tilde K_tW_{} \rangle\\
&= \langle W_{x},\tilde KW_{t} \rangle +\frac{1}{2}\langle
W_{},\big(\tilde K_x+ (\alpha_x/\alpha)\big) W_{t} \rangle
+\frac{1}{2}\langle W_{x},\tilde K_tW_{} \rangle ,\\
\end{aligned}
\label{parts} \end{equation}
 we obtain by calculations similar to the above the
auxiliary energy estimate:
\begin{equation}
\begin{aligned} \frac{1}{2}\langle W_{x},\tilde KW_{} \rangle_t &\le
-\langle  W_{x}, \tilde K(\tilde A^0)^{-1}\tilde A W_{x}\rangle\\
&\qquad +C(C_*)|w^{II}_x|^2
+C\langle (|\bar W_x|+ \bar \zeta + \zeta) |w^I_x|,|w^I_x|\rangle\\
&\qquad+ C\bar \zeta^{-1}|w^{II}_{xx}|^2
+ C(C_*)(|W|_{L^2}+|\xi|^2_{H^1}) ,
\end{aligned}
\label{est2} \end{equation}

 where $\bar \zeta>0$ is an arbitrary constant
arising through Young's inequality. (Here, we have estimated term
$\langle \tilde A U_x,(\alpha_x/\alpha)U\rangle$ arising in the
middle term of the righthand side of (\ref{parts}) using
(\ref{symm1}) by $C(C_*)\int |\bar W_x||U|^2\le
C(C_*)|U|_{L^\infty}^2$.)

\smallskip
{\bf Combined, weighted $H^1$ estimate.} Choosing $ \zeta << \bar
\zeta << 1$, adding (\ref{est2}) to the sum of (\ref{nextsym})
times a suitably large positive constant $M(C_*,\bar \zeta)>> \bar
\zeta^{-1}$, and (\ref{sym}) times  $M(C_*,\bar \zeta)^2$ and
recalling \ref{skew}, we obtain, finally, the combined first-order
estimate
\begin{equation}
\begin{aligned} \frac{1}{2} \Big( &M(C_*,\bar \zeta)^2\langle W, \tilde A^0
W_{}\rangle + \langle W_{x},\tilde KW_{} \rangle
+M(C_*,\bar \zeta)\langle W_x, \tilde A^0 W_{x}\rangle \Big)_t \\
&\le -\theta (|W_{x}|_{L^2}^2 +|w^{II}_{xx}|_{L^2}^2) +C(C_*)
\left(|W|_{L^2}^2+|\xi|^2_{H^1}
 \right),\\
\end{aligned}
\label{estfinal} \end{equation}

 $\theta>0$, for any $\bar \zeta$, $\zeta(\bar
\zeta,C_*)$ sufficiently small, and $C_*$, $C(C_*)$ sufficiently
large.

\medskip
{\bf Higher order estimates.} Performing the same procedure on the
twice- and thrice-differentiated versions of equation
(\ref{leibnitz}), we obtain, likewise, Friedrichs estimates
\begin{equation}
\begin{aligned} \frac{1}{2}\langle \partial_x^qW_,\tilde A^0 &\partial_x^q
W\rangle_t \le -(\theta/2)|\partial_x^{q+1} w^{II}|_{L^2}^2 -
(C_*\theta/4) \langle |\bar W_x||\partial_x^qw^I|,|\partial_x^qw^I|\rangle\\
& + C(C_*)\big( \zeta |\partial_x^qw^I|_{L^2}^2 +
|\partial_x^qw^{II}|_{L^2}^2 + |W_x|_{H^{q-2}_\alpha}
+|W|_{L^2}^2+ |\xi|^2_{H^{2q}}\big),\\
\end{aligned}
\label{symm23}
\end{equation}
 and Kawashima estimates
\begin{equation}
\begin{aligned} \frac{1}{2}\langle \partial_x^q W, \tilde
K&\partial_x^{q-1} W \rangle_t \le
-\langle  \partial_x^q W, \tilde K(\tilde A^0)^{-1}\tilde A \partial_x^q W\rangle\\
&+C(C_*)|\partial_x^q w^{II}|_{L^2}^2
+C\langle (|\bar W_x|+ \bar \zeta + \zeta) |\partial_x^q w^I|,|\partial_x^q w^I|\rangle\\
&+ C\bar \zeta^{-1}|\partial_x^{q+1} w^{II}|_{L^2}^2 +
C(C_*)(|W_x|_{H^{q-2}_\alpha} + |W|_{L^2}^2+ |\xi|^2_{H^{2q}}),
\end{aligned}
\label{est23} \end{equation}
 for $q=2, \, 3$, provided $\bar \zeta$, $\zeta(\bar
\zeta,C_*)$ are sufficiently small, and $C_*$, $C(C_*)$ are
sufficiently large. The calculations are similar to those carried
out already; see also the closely related calculations of Appendix
A, \textbf{\cite{MaZ.2}.}

\smallskip
{\bf Final estimate.} Adding $M(C_*,\bar \zeta)^2$ times
(\ref{estfinal}), $M(C_*,\bar \zeta)$ times (\ref{symm23}), and
(\ref{est23}), with $q=2$, where $M$ is chosen still larger if
necessary, we obtain
\begin{equation}
\begin{aligned} \frac{1}{2} \Big( &M(C_*,\bar \zeta)^4\langle W, \tilde
A^0 W_{}\rangle + M(C_*,\bar \zeta)^2\langle W_{x},\tilde KW_{}
\rangle
+M(C_*,\bar \zeta)^3\langle W_x, \tilde A^0 W_{x}\rangle \\
&\qquad + \langle \partial_x^2W,\tilde K\partial_xW \rangle
+M(C_*,\bar \zeta)\langle \partial_x^2W,
\tilde A^0 \partial_x^2W\rangle \Big)_t \\
&\le -\theta (| W_x|_{H^1_\alpha}^2 +|w^{II}_x|_{H^2_\alpha}^2)
+C(C_*) \left(|W|_{L^2}^2) + |\xi|^2_{H^4} \right).\\
\end{aligned}
\label{estfinalq2} \end{equation}

Adding now $M(C_*,\bar \zeta)^2$ times (\ref{estfinalq2}),
$M(C_*,\bar \zeta)$ times (\ref{symm23}), and (\ref{est23}), with
$q=3$, we obtain the final higher-order estimate
\begin{equation}
\begin{aligned} \frac{1}{2} \Big( &M(C_*,\bar \zeta)^6\langle W, \tilde
A^0 W_{}\rangle + M(C_*,\bar \zeta)^4\langle W_{x},\tilde KW_{}
\rangle
+M(C_*,\bar \zeta)^5\langle W_x, \tilde A^0 W_{x}\rangle \\
&\qquad + M(C_*,\bar \zeta)^2 \langle \partial_x^2W,\tilde
K\partial_xW \rangle +M(C_*,\bar \zeta)^3\langle \partial_x^2W,
\tilde A^0 \partial_x^2W\rangle \\
&\qquad + \langle \partial_x^3W,\tilde K\partial_x^2W \rangle
+M(C_*,\bar \zeta)\langle \partial_x^3W,
\tilde A^0 \partial_x^3W\rangle \Big)_t \\
&\le -\theta (| W_x|_{H^2_\alpha}^2 +|w^{II}_x|_{H^2_\alpha}^2)
+C(C_*) \left(|W|_{L^2}^2 + |\xi|^2_{H^6} \right).\\
&\le -\theta | W|_{H^3_\alpha}^2
+C(C_*) \left(|W|_{L^2}^2) + |\xi|^2_{H^6} \right).\\
\end{aligned}
\label{estfinalq2} \end{equation}

Denoting

\begin{equation}
\begin{aligned} \mathcal{E} (W):= \frac{1}{2} \Big( M(C_*,\bar \zeta)^6\langle
W, \tilde A^0 W_{}\rangle + M(C_*,\bar \zeta)^4\langle
W_{x},\tilde KW_{} \rangle\\
 +M(C_*,\bar \zeta)^5\langle W_x, \tilde A^0 W_{x}\rangle
 + M(C_*,\bar \zeta)^2 \langle \partial_x^2W,\tilde
K\partial_xW \rangle +\\
M(C_*,\bar \zeta)^3\langle \partial_x^2W, \tilde A^0
\partial_x^2W\rangle + \langle \partial_x^3W,\tilde
K\partial_x^2W \rangle +M(C_*,\bar \zeta)\langle \partial_x^3W,
\tilde A^0 \partial_x^3W\rangle \Big), \\
\end{aligned}
\end{equation}

 we have by Young's inequality that $\mathcal{E}^{1/2}$ is equivalent
to norms $H^3$ and $H^3_\alpha$, hence (\ref{estfinal}) yields
$$
\mathcal{E}_t\le -\theta_2 \mathcal{E} + C(C_*) \left(|W|_{L^2}^2)
+ |\xi|^2_{H^6} \right),
$$
from which we conclude,
$$
\mathcal{E}(t)\le e^{-\theta_2 t}\mathcal{E}(0) + C(C_*) \int_0^t
e^{-\theta_2 (t-s)} \left(|W|_{L^2}^2) + |\xi|^2_{H^6} \right)
(s)\, ds.
$$
This is equivalent to (\ref{whbounds}).

{\bf The general case.} It remains only to discuss the general
case that hypotheses (A1)--(A3) hold as stated and not everywhere
along the profile, with $\tilde G$ possibly nonzero. These
generalizations requires only a few simple observations. The first
is that we may express matrix $\tilde A$ in (\ref{leibnitz}) as
$$
\tilde A=
{\hat A} + (|\bw_x|+ \zeta)\left( \begin{matrix} 0&\CalO(1)\\
\CalO(1) & \CalO(1)\\  \end{matrix} \right), \label{symmed}
$$
where $\hat A$ is a symmetric matrix obeying the same derivative
bounds as described for $\tilde A$, identical to $\tilde A$ in the
$11$ block and obtained in other blocks $jk$ by smoothly
interpolating over a bounded interval $[-R,+R]$ between $\bar
A(\tilde W_{-\infty})_{jk}$ and $\bar A(\tilde W_{+\infty})_{jk}$.
Replacing $\tilde A$ by $\hat A$ in the $q$th order
Friedrichs-type bounds above, we find that the resulting error
terms may be expressed as (integrating by parts if necessary)
$$
\langle \partial_x^q \CalO(|\bw_x|+\zeta) |W|,
|\partial_x^{q+1}w^{II}| \rangle
$$
plus lower-order terms, hence absorbed using Young's inequality to
recover the same Friedrichs-type estimates obtained in the
previous case. Thus, we may relax (A1') to (A1).

The second observation is that, because of the favorable terms
$$
-(C_*\theta/4) \langle
|\bw_x||\partial_x^qw^I|,|\partial_x^qw^I|\rangle
$$
occurring in the righthand sides of the Friedrichs-type estimates,
we need the Kawashima-type bound only to control the contribution
to $|\partial_x^q w^I|^2$ coming from $x$ near $\pm \infty$; more
precisely, we require from this estimate only a favorable term
$$
-\theta \langle \big(1- \CalO(|\bw_x|+\zeta + \bar \zeta)\big)
|\partial_x^qw^I|,|\partial_x^qw^I|\rangle \label{weakKaw}
$$
rather than $ -\theta |\partial_x^qw^I|^2$ as in (\ref{est2}) and
(\ref{est23}). But, this may easily be obtained by substituting
for $\tilde K$ a skew-symmetric matrix-valued function $\hat K$
defined to be identically equal to $\bar K(+\infty)$ and $\bar
K(-\infty)$ for $|x|>R$, and smoothly interpolating between $\bar
K(\pm\infty)$ on $[-R,+R]$, and using the fact that
$$
\big(\bar K (\bar A^0)^{-1}\bar A + \bar B\big)_\pm \ge \theta >
0,
$$
hence
$$
\big(\hat K (\tilde A^0)^{-1}\tilde A + \tilde B\big) \ge \theta
(1- \CalO(|\bw_x|+\zeta)).
$$
Thus, we may relax (A2') to (A2).

Finally, notice that the term $\tilde G-\bar G$ in the
perturbation equation may be Taylor expanded as
$$
\left( \begin{matrix} 0 \\ \tilde g(\tilde W_x, \bar U_x) + g(\bw_x, \tilde W_x)
 \end{matrix} \right) + \left( \begin{matrix} 0 \\
\CalO(|W_x|^2) \end{matrix} \right) \label{gexp}
$$
The first, linear term on the righthand side may be grouped with
term $\tilde A^0 W_x$ and treated in the same way, since it decays
at plus and minus spatial infinity and vanishes in the $1$-$1$
block. The $(0,\CalO(|W_x|^2)$ nonlinear term may be treated as
other source terms in the energy estimates Specifically, the
worst-case terms $\langle \partial_x^3 W, K \partial_x^2
\CalO(|W_x|^2)\rangle$ and $ \langle \partial_x^3 W,  \partial_x^3
(0, \CalO(|W_x|^2))\rangle = \langle \partial_x^4 w^{II},
\partial_x^2 \CalO(|W_x|^2)\rangle $ may be bounded, respectively,
by $|W|_{W^{2,\infty}}|W|_{H^3}^2$ and
$|W|_{W^{2,\infty}}|w^{II}|_{H^4}|W|_{H^3}$. Thus, we may relax
(A3') to (A3), completing the proof of the general case (A1)--(A3)
and the theorem.
\end{proof}

\end{document}